\documentclass[10pt, reqno]{amsart}

\usepackage[nobysame,abbrev,alphabetic]{amsrefs}

\usepackage{amsmath, amsthm, amsfonts, amssymb, latexsym}
\usepackage[all]{xy}
\usepackage[hidelinks]{hyperref}
\usepackage{cleveref}
\usepackage{microtype}
\usepackage{verbatim}
\usepackage[normalem]{ulem}

\usepackage{tikz}
\usetikzlibrary{matrix, arrows, decorations.pathmorphing}
\usepackage{tikz-cd}
\usetikzlibrary{arrows}

\usepackage{mathrsfs}

\usepackage{epsfig}
\usepackage{amscd}
\usepackage{enumerate}

\newcommand{\Proof}[1]{\begin{proof} #1 \end{proof}}
\newtheorem{thm}{Theorem}[section]

\newtheorem{pro}[thm]{Proposition}
\newtheorem*{thmm}{Theorem}

\newtheorem*{proo}{Proposition}

\newcommand{\Pro}[1]{\begin{pro} #1\end{pro}}
\newtheorem{cor}[thm]{Corollary}
\newcommand{\Cor}[1]{\begin{cor} #1 \end{cor}}
\newtheorem{lem}[thm]{Lemma}
\newcommand{\Lem}[1]{\begin{lem} #1 \end{lem}}
\newtheorem{rem}[thm]{Remark}

\theoremstyle{definition}
\newtheorem{defn}[thm]{Definition}

\newtheorem{Ex}[thm]{Example}

\newtheorem{proposition}[equation]{Proposition}

\theoremstyle{remark} 
\newtheorem{remark}[equation]{Remark}

\newcommand{\bF}{{\mathbb{F}}}
\newcommand{\bZ}{{\mathbb{Z}}}
\newcommand{\Inf}{\operatornamewithlimits{inf}}
\newcommand{\Sq}{\operatorname{\sf Sq}}

\renewcommand{\ge}{\geqslant}

\usepackage[parfill]{parskip}
\usepackage[margin=1in]{geometry}
\numberwithin{equation}{section}

\newcommand{\rar}{\ensuremath{\rightarrow}}
\newcommand{\lrar}{\ensuremath{\longrightarrow}}

\newcommand{\la}{\langle}
\newcommand{\ra}{\rangle}

\newcommand{\stk}[1]{\stackrel{#1}{\rightarrow}}
\newcommand{\lstk}[1]{\stackrel{#1}{\longrightarrow}}

\newcommand{\Gal}{\textup{Gal}}

\renewcommand{\inf}{\textup{inf\,}}
\newcommand{\fp}{\mathbb{F}_p}
\newcommand{\fq}{\mathbb{F}_q}

\newcommand{\Dec}{\textup{Dec}}
\newcommand{\sep}[1]{{#1}_{\textup{sep}}}
\newcommand{\ints}{\mathbb{Z}}
\newcommand{\Hom}{\text{Hom}}
\newcommand{\ZZ}{\mathbb{Z}}
\newcommand{\FF}{\mathbb{F}}
\newcommand{\Ss}{\mathbb{S}}
\newcommand{\Span}[1]{\ensuremath{ \langle #1 \rangle }}
\newcommand{\res}{\text{res}}
\newcommand{\set}[1]{\ensuremath{ \lbrace #1 \rbrace }}
\DeclareMathOperator{\Corres}{cor}

\newcommand{\dirlim}[1]{\ensuremath{\lim_{\substack{\longrightarrow \\
{\scriptscriptstyle #1}}}}\;}

\newcommand{\invlim}[1]{\ensuremath{\lim_{\substack{\longleftarrow \\
{\scriptscriptstyle #1}}}}\;}

\newcommand{\dlim}[1]{\ensuremath{\underset{\underset{#1}{\rightarrow}}{\text{lim}}\;}}
\newcommand{\ilim}[2]{\ensuremath{\underset{\underset{#1}{\leftarrow}}{\text{lim}^{#2}}\;}}

\makeatletter
\newcommand{\addresseshere}{
  \enddoc@text\let\enddoc@text\relax
}
\makeatother

\begin{document}

\title{The Bloch--Kato conjecture, decomposing fields, and generating cohomology in degree one}
\date{\today}

\author{Sunil K. Chebolu}
\address[Chebolu]{Department of Mathematics \\
Illinois State University\\
Normal, IL 61761 USA}
\email{schebol@ilstu.edu}

\author{J\'{a}n Min\'{a}\v{c}}
\address[Min\'{a}\v{c}]{Department of Mathematics\\
The University of Western Ontario\\
London, ON N6A 5B7, Canada}
\email{minac@uwo.ca}

\author{C\.{i}han Okay}
\address[Okay]{Department of Mathematics, Bilkent University, 06800 Ankara, Turkey}
\email{cihan.okay@bilkent.edu.tr}

\author{Andrew Schultz}
\address[Schultz]{Wellesley College, Wellesley, MA 02481 USA}
\email{andrew.c.schultz@gmail.com}

\author{Charlotte Ure}
\address[Ure]{Department of Mathematics, Illinois State University,
Normal, IL 61761 USA}
\email{cure@ilstu.edu}

\makeatletter
\let\@wraptoccontribs\wraptoccontribs
\makeatother
\contrib[with an appendix by]{David Benson}

\thanks{Sunil Chebolu is partially supported by Simons Foundation’s Collaboration Grant for Mathematicians (516354). J\'{a}n Min\'{a}\v{c} is partly supported by the National Sciences and Engineering Research Council of Canada (NSERC) grant R0370A01. He also gratefully acknowledges the Western University Faculty of Science Distinguished Professorship for 2020-2021 and the support of the Western Academy for Advanced Research, where he was elected as a Fellow for the 2022-2023 program. Charlotte Ure was partially supported by an AMS--Simons Travel Grant and by a New Faculty Initiative Grant from Illinois State University. }



\dedicatory{Dedicated to the memory of Vladimir Voevodsky with admiration and gratitude.}

\begin{abstract}
The famous Bloch--Kato conjecture implies that for a field $F$ containing a primitive $p$th root of unity, the cohomology ring of the absolute Galois group $G_F$ of $F$ with $\fp$ coefficients is generated by degree one elements. We investigate other groups with this property and characterize all such groups that are finite. Restricting to the case 
of $p$-groups, our work answers a question of Quadrelli, Snopce and Vanacci posed in 2022. As a further step in this program, we study implications of the Bloch--Kato conjecture to cohomological
invariants of finite field extensions. Conversely, these cohomological invariants have implications for refining the Bloch--Kato conjecture. In service of such a refinement, we define the notion of a decomposing field for a cohomology class of a finite field extension and study minimal decomposing fields of degree two cohomology classes arising from degree $p$ extensions. We illustrate this refinement by explicitly computing the cohomology rings of superpythagorean fields and $p$-rigid fields. Finally, we construct nontrivial examples of cohomology classes and their decomposing fields, which rely on computations by David Benson in the appendix.  
\end{abstract}

\maketitle


\section{Introduction} 
In this paper, we present a new approach to refining the celebrated Bloch--Kato conjecture and also obtain this refinement for cohomology classes of degree two. We begin by recalling the Bloch--Kato conjecture. 

Let $F$ be a field that contains a primitive $p$-th root of unity, where $p$ is a fixed prime number. The absolute Galois group of $F$ is
the Galois group $G_{F}:= \Gal(F_{sep}/F)$,
where $F_{sep}$ is the separable closure of $F$. This is a rich profinite group whose structure is mysterious. To better understand this profinite group, one
studies its continuous, or Galois, cohomology $H^{*}(G_{F}, \fp)$, where we denote by $\fp$ the trivial $G_F$-module with $p$ elements. This is a graded commutative ring,
and the Bloch--Kato conjecture describes this ring in terms of generators and relations. In particular, it is generated in degree one and all its relations are generated in degree two. A graded ring which is generated in degree one and whose relations are generated in degree two is called a \emph{quadratic algebra}, and hence the Bloch--Kato conjecture tells us that the cohomology ring $H^*(G_F,\mathbb{F}_p)$ is a quadratic algebra. This is an extraordinary statement about the structure of the Galois cohomology ring. 

This property helps distinguish those profinite groups that arise as the absolute Galois groups of a field from those profinite groups that do not. The characterization of absolute Galois groups among profinite groups is one of the major open questions in the current study of Galois theory. It is very interesting to determine those profinite groups $G$ for which $H^*(G, \fp)$ is a quadratic algebra. For a recent paper addressing this question, see \cite{Q22}. Our approach to this question is complementary to the approach in \cite{Q22}. We view our profinite group as a projective limit of other profinite groups, and we consider $H^* (G,\fp)$ as an inductive limit of corresponding cohomology rings.
In the case when $H^* (G,\fp)$ is a quadratic algebra, it is often a thrilling detective mystery to clarify precisely when certain elements become fully decomposable under inflation in the process of passing to the inductive limit. 
(See \Cref{defn:fully-decomposes} and \Cref{defn:decomposing-field} below.)

As progress towards this larger goal, in the theorem below we characterize all \emph{finite} groups $G$ which satisfy the following slightly weaker condition:
\begin{equation}
H^*(G,\mathbb{F}_p) \text{ is generated in degree one.}
\tag{$\star$}
\end{equation}
This  answers a question posed by Quadrelli, Snopce, and Vannacci \cite[Remark 2.10]{Q22} (see also \Cref{Claudio}).
\begin{thm}\label{thm:intro:cohomology.generated.in.degree.one}
    Suppose that $G$ is a finite group, $p$ is a prime number that divides $|G|$, and $G$ acts trivially on $\fp$. Then $H^*(G,\fp)$ is generated in degree one if and only if $p=2$ and the Sylow $2$-subgroup of $G$ is a nontrivial elementary $2$-abelian group that admits a normal complement.
\end{thm}
Note that when $n$ is an odd number, the dihedral group of order $2n$ is an example of a group which satisfies the conditions of this theorem. It is well-known that the only non-trivial finite group which arises as the absolute Galois group of a field is the cyclic group of order $2$, and hence the study of those groups satisfying condition ($\star$) provides an interesting liminal property in the investigation of absolute Galois groups. For more equivalent conditions and a proof of this theorem, see \Cref{thm:cohomology.generated.degree.one}. 

We propose a systematic study of pro-$p$ groups $G$ that satisfy ($\star$) or the stronger property that $H^*(G,\mathbb{F}_p)$ is a quadratic algebra using methods in group cohomology of finite groups and studies of Galois embedding problems. This new point of view makes our paper nicely complementary with \cite{CEM} and \cite{Q22} where some similar problems are investigated using different tools. 

In the case of $H^*(G_F,\mathbb{F}_p)$, the profinite structure of $G_F$ suggests some deep questions about precisely how any given $\alpha \in H^*(G_F,\mathbb{F}_p)$ is generated from $H^1(G_F,\mathbb{F}_p)$.  Specifically, because $G_F$ is profinite, we know that there must exist some finite Galois extension $K/F$ and some $\gamma \in H^*(\text{Gal}(K/F),\mathbb{F}_p)$ with $\alpha = \inf_{\text{Gal}(K/F)}^{G_F}(\gamma)$.  On the other hand, there is no reason to expect that $\gamma$ is in the subalgebra of $H^*(\text{Gal}(K/F),\mathbb{F}_p)$ generated by $H^1(\text{Gal}(K/F),\mathbb{F}_p)$.  (Indeed, Theorem \ref{thm:intro:cohomology.generated.in.degree.one} tells us $H^*(\text{Gal}(K/F),\mathbb{F}_p)$ will contain elements outside the subalgebra generated in degree one ``most of the time".) It is natural to wonder, then, precisely when along the tower $F_{sep}/F$ the inflation of $\gamma$ takes this form.  This motivates the following definition.

\begin{defn}\label{defn:fully-decomposes}
Suppose that $K$ is a finite Galois extension of $F$ and  $\gamma \in H^m(\Gal(K/F),\fp)$ for some $m \geq 2$. Then $\gamma$ \emph{is fully decomposable} if there exists some $\ell \in \mathbb{N}$ and some elements $\{\beta_{i,j}: 1 \leq i \leq \ell, 1 \leq j \leq m\} \subseteq H^1(\Gal(K/F),\fp)$ so that 
\begin{equation}\label{eqn:decompose}
    \gamma = \sum_{i=1}^\ell \beta_{i,1} \cup \cdots \cup \beta_{i,m}.
\end{equation} 
Otherwise, we say that $\gamma$ is \emph{not fully decomposable.} 
Furthermore, if $L/F$ is a Galois extension with $K \subseteq L$, then we say that $\gamma$ \emph{fully decomposes over $L$} if  $\inf_{\Gal(K/F)}^{\Gal(L/F)}(\gamma)$ is fully decomposable in $H^m( \Gal(L/F), \fp)$. Otherwise, we say that $\gamma$ is \emph{not fully decomposable over $L$}.
\end{defn}

We remark that if $\gamma$ becomes trivial over $L$, then it also fully decomposes over $L$. We construct an example of finite field extensions $F\subset K \subset L$  and a class $\gamma \in H^4\left( \Gal(K/F), \FF_2\right)$ so that $\gamma$ fully decomposes over $L$ --- but $\inf_{\Gal(K/F)}^{\Gal(L/F)}(\gamma) \neq 0$  --- in \Cref{sec:splitting.without.vanishing}. This example relies on computations by David Benson in \Cref{appendix} for the inflation of a nontrivial cohomology class in $H^4\left( Q_8, \FF_2\right)$, where $Q_8$ is the quaternion group with eight elements. 

To study the decomposition of some fixed $\gamma \in H^m\left( \Gal(K/F), \FF_p\right)$ systematically, we propose the following definition. 

\begin{defn}[Decomposing field]\label{defn:decomposing-field}
Suppose that $K$ is a finite Galois extension of $F$ and that $\gamma \in H^m(\Gal(K/F),\fp)$ for some $m \geq 2$. A \emph{decomposing field of $\gamma$} is a field over which $\gamma$ fully decomposes. We call it a \emph{minimal decomposing field} if it is a decomposing field of minimal degree. 
\end{defn} 

One key observation (see \Cref{prop:finite-Bloch--Kato}) to make is that the following statements are equivalent: (1) for every finite Galois extension $K/F$ and each $\gamma \in H^m(\Gal(K/F),\fp)$ there exists a decomposing field for $\gamma$ that is finite degree over $F$, and (2) the Galois cohomology ring $H^m(G_F,\fp)$ is generated by degree one elements. For this reason, we propose the study of decomposing fields as a method to provide a more refined understanding of the Bloch--Kato conjecture.

In addition to giving deeper insight into the heartbeat of the Bloch--Kato conjecture, decomposing fields suggest additional interesting questions. For example, the degree of a minimal decomposing field of $\gamma$ is an interesting numeric invariant to study.  Another such invariant would be the minimal number $\ell$ for which we can find $\{\beta_{i,j}: 1 \leq i \leq \ell, 1 \leq j \leq m \} \subseteq H^1(\Gal(L/F),\fp)$ satisfying \Cref{eqn:decompose}. This value $\ell$ is known as the \emph{symbol length} of the element. Those $\gamma$ that admit such an expression when $\ell=1$ can be considered as ``pure decompositions" as a product of $1$-cohomology cycles. These decompositions are of special interest.

In this paper, we begin an investigation into decomposing fields by studying the Galois tower of fields of $F_{sep}$ given by
\begin{equation}\label{eq:filtration.tower} F = F^{(1)} \subset F^{(2)} \subset F^{(3)} \subset \cdots \subset F^{(n)} \subset \cdots \subset F(p) \subset F_{sep},
\end{equation}
where $F^{(2)}$ is the compositum of all cyclic extensions of degree $p$ over $F$; for $n \ge 2$,
$F^{(n+1)}$ is the compositum of all cyclic extensions of degree $p$ over $F^{(n)}$ which
are Galois over $F$; and $F(p)$ is the maximal $p$-extension of $F$. Using basic properties of $p$-groups and Galois theory, it can be shown that $F(p) = \bigcup_{n\geq 1} F^{(n)}$. We will call the filtration from (\ref{eq:filtration.tower}) ``the filtration tower associated with $F$.'' We remark that as a consequence of the Bloch--Kato conjecture, the inflation $H^*(\Gal(F(p)/F), \fp ) \rightarrow H^*(G_F, \fp)$ is an isomorphism (see \Cref{pro:Cohology-F(p)}). 

For each positive integer $n$, we denote the Galois group of $F^{(n)}$ over $F$ by $G_F^{[n]}$. We remark that $G_F^{[n]}$ can also be thought of as $G_F/G_F^{(n)}$, where $G_F^{(n)}$ is the $n$th term in the lower $p$-central series of $G_F$.

As a first result in the study of decomposing fields, we investigate the decomposition of cohomology classes in $H^*( G_F^{[n]}, \fp)$. In \cite{CEM}, it was shown that any cohomology class from $H^2(G_F^{[n]},\fp)$ fully decomposes over $F^{(n+1)}$. We first reprove this statement using the tools of Brauer embedding theory in \Cref{quotients} and then prove the following refinement of it.  

\begin{thm}\label{thm:main}
Let $F$ be a field that contains a primitive $p$-th root of unity, and let $n \ge 3$. Let $\gamma \in H^2(G_F^{[n]},\fp)$ be an element that is not fully decomposable. Then a minimal decomposing field for $\gamma$ is a Galois extension of degree $p$ over $F^{(n)}$ that is also Galois over $F$.\end{thm} 

In addition to this refinement, we give the following alternate characterization of $F^{(n+1)}$ using the notion of decomposing fields. 

\begin{thm}\label{thm:theorem.after.main}
Let $n \geq 3$. The field $F^{(n+1)}$ is a field minimal with respect to inclusion that is a decomposing field for all $\gamma \in H^2(G_F^{[n]},\fp)$.
\end{thm}

These two theorems are proven in Section \ref{refinement}, with the previous theorem made more precise in \Cref{thm:F^(n+1)minimal-decomposing}.

We showcase Theorem \ref{thm:main} in the case of superpythagorean fields and $p$-rigid fields for $p$ an odd prime. For precise definitions, see \Cref{computation}. We shall first describe the case where $G_F(p)=\Gal(F(p)/F)$ is finitely generated, both in the case where $F$ is superpythagorean (when we set $p=2$) and when $F$ is $p$-rigid (in which case we assume $p$ is odd).  For these classes of fields and under this assumption of finite generation, the Galois group $G_F(p)$ is the inverse limit of the semidirect products of certain abelian $p$-groups. Specifically, define 
\begin{align*}
A_2(n;d) &=  (C_{2^n})^d\rtimes C_2, &&\text{ with } xyx^{-1} = y^{-1} \text{ for } y \in (C_{2^n})^d, \langle x\rangle = C_2; \text{ and}\\
B_p(n;d,k) &= (C_{p^n})^d\rtimes C_{p^n}, &&\text{ with } xyx^{-1}=y^{p^k+1} \text{ for }y \in (C_{p^n})^d  \text{ and }   x \text{ a fixed generator of }C_{p^n},
\end{align*}
where for a given $\ell \in \mathbb{N}$ we have used the notation $C_\ell$ to denote the cyclic group of order $\ell$.  When $F$ is superpythagorean, then there exists $d \in \mathbb{N}$ so that if we define $P(n)=A_2(n;d)$ then it follows that $G_F(2) = \varprojlim P(n)$.  Likewise if $F$ is $p$-rigid for some odd prime $p$, then there exist $d \in \mathbb{N}$ and $k \in \mathbb{N} \cup\{\infty\}$ so that if we define $P(n)=B_p(n;d,k)$ then it follows that $G_F(p)= \varprojlim P(n)$. (By convention, in the case $k=\infty$ the relation $xyx^{-1}=y^{p^k+1}$ is interpreted as $xyx^{-1}=y$.)

The case where $F$ is superpythagorean or $p$-rigid but with $G_F(p)$ infinitely generated has a natural connection to the finitely generated case, since in both cases $G_F(p)$ can be expressed as the projective limit of groups of the types described above. Specifically, if we define $P(n)=\ilim{d}{} A_2(n;d)$ in the superpythagorean case and $P(n)=\ilim{d}{} B_p(n;d,k)$ in the $p$-rigid case, then we again recover $G_F(p)$ as $\ilim{}{} P(n)$. Here, the projective limit $\ilim{d}{}A_2(n;d)$ represents the limit over all quotients of $G_F(2)$ of the form $A_2(n;d)$ for some $d \in \mathbb{N}$. The limit $\ilim{d}{} B_p(n;d,k)$ should be interpreted in a similar way.

The main goal is to examine the inflation map arising from this inverse system of groups $P(n)$, and explicitly verify --- using group cohomological methods --- the following instance of \Cref{thm:main}. 
\begin{thm}\label{thm:example}
    Let $F$ either be a superpythagorean field, or a $p$-rigid field containing a primitive $p$-th root of unity  for $p$ an odd prime. Furthermore, let $P(n)$ be defined as above.  Then every $\gamma \in H^2(P(n),\fp)$ fully decomposes in $H^2(P(n+1),\fp)$ under the inflation map.  
\end{thm}

Whereas the proof of Theorem \ref{thm:main} uses the Bloch--Kato conjecture, our proof of Theorem \ref{thm:example} is completed independently of this powerful result. Indeed, the proof relies on the calculation of Galois groups for superpythagorean and $p$-rigid fields (see \cite{OtherBecker}, \cite[Ch.~3]{Becker}, \cite[Thm.~4.10]{CMQ}, \cite[Cor.~4.10]{MRT16}, \cite{WareII}) and, subsequently, the computation of group cohomology for these groups (see section \ref{subsec:cohomology.of.metacyclic.groups}), neither of which uses the Bloch--Kato conjecture. 

\noindent
\textbf{Outline of paper.}
Our paper is organized as follows. In \Cref{tools}, we recall essential tools which will be used in our work, such as the Bloch--Kato conjecture, Kummer theory, and Galois embeddings. \Cref{quotients} revisits the main theorem of \cite{CEM}, offering a different proof and setting the stage for the later results. The main results on minimal decomposing fields and a refinement of the Bloch--Kato conjecture (Theorems \ref{thm:main} and \ref{thm:theorem.after.main}) are proved in \Cref{refinement}. We illustrate our results with examples in Section \ref{sec:splitting.without.vanishing} and \ref{computation}. The former, together with \Cref{appendix}, offers explicit examples of quotients of finite groups and cohomology classes that decompose fully but do not vanish under the inflation map. \Cref{computation} highlights the refinement of the Bloch--Kato conjecture in the example of superpythagorean and $p$-rigid fields. 
Finally, we characterize finite groups whose cohomology ring is generated in degree $1$ in \Cref{sec:BKpropgroups}.

\noindent
\textbf{Acknowledgments.} 
We thank John Swallow for sharing his insights in our early stage of this project. 
We also extend our gratitude to Stefan Gille, Mathieu Florence, and Johannes Huebschmann for their insightful discussions, and to Ido Efrat and Claudio Quadrelli for their collaborative work, all of which inspired some of the ideas presented in this paper. We are very grateful to David Benson for a number of comments, suggestions, and group cohomological calculations. He informed us of a result of Swan that helped resolve a question in group cohomology; see \Cref{thm:cohomology.generated.degree.one}. We are also very grateful to David Benson for writing an appendix to this paper. Finally, we would like to thank two anonymous referees of this paper for their careful review of our initial manuscript and for providing helpful comments. 

\section{Preliminaries} \label{tools}
For basic facts about Galois cohomology (which is just the continuous group cohomology of the profinite group $G_{F}$), we refer 
the reader to \cites{NeukirchSchmidtWingberg, Serre65, SerreCG}. In this section, we explain some of the arithmetic tools and techniques used in the remainder of the paper. Our toolkit consists of the Bloch--Kato conjecture, Kummer theory, Galois embedding theory, and the fundamental five-term sequence. 

\subsection{History and basic facts on the Bloch--Kato conjecture}

To explain the precise statement of the Bloch--Kato conjecture, we have to introduce another ring called reduced Milnor $K$-theory, whose structure is relatively more transparent. For a field $F$, we write $F^\times$ for the multiplicative group of non-zero elements of $F$. 
The reduced Milnor $K$-theory, denoted $k_{*}F$, is  the quotient of the tensor algebra $T(F^\times)$ by the two-sided ideal  $\la a \otimes b \, | \, a + b = 1,\ \  a, b  \in F^{\times} \ra$ reduced modulo $p$. That is,
\begin{equation}\label{eqn:k_*F} k_{*}F :=  \frac{T(F^{\times})}{ \la a \otimes b \, | \, a + b = 1, \ \  a, b  \in F^{\times} \ra} \otimes \fp. \end{equation}
This is a graded ring by setting the degree of any $x \in F^\times/(F^\times)^p$ to $1$. Using Kummer theory and a lemma of Bass and Tate \cite{Milnor70}*{Lemma 6.1}, it is possible to construct (see next section) a natural graded-algebra homomorphism $\eta$ called the norm-residue map from reduced Milnor $K$-theory to Galois cohomology:
\begin{align}\label{eq:norm.residue.homomorphism} 
\eta \colon k_{*}F \longrightarrow H^{*}(G_{F}, \fp).
\end{align}
The Bloch--Kato conjecture is the statement that this map is an isomorphism of rings. It is known that this map is well-defined, and we will reprove it for the reader's convenience using Galois embedding problems in the proof of \Cref{thm:G_F^{[3]}}. It has been shown \cite[Page 169]{Voe-Sus} that the surjectivity of this map implies injectivity. Note that since elements of $F^{\times}/(F^{\times})^p$ have degree one in $k_*(F)$, the Bloch--Kato conjecture implies the aforementioned statements about the generators being in degree one and relations in degree two, respectively.

The Bloch--Kato conjecture was resolved in 2011 by Voevodsky and Rost, but there is a long and interesting history leading to its full proof. 
Merkurjev and Suslin showed in \cite{MerkurjevSuslin82} that (\ref{eq:norm.residue.homomorphism}) is an isomorphism in degree two, in which case the right-hand side may be identified with the Brauer group. The case $p=2$ was implicitly suggested by Milnor in \cite{Milnor70} in 1970, and it was eventually proved by Voevodsky in \cite{Voevodsky03a}. Voevodsky and Rost later proved the case when $p$ is odd; for this reason, the Bloch--Kato conjecture is often referred to as the Rost--Voevodsky theorem. 
In this paper we will often use the name ``Bloch--Kato conjecture" (or simply ``Bloch--Kato") when referring to the fact that $\eta$ from Equation (\ref{eq:norm.residue.homomorphism}) is an isomorphism of rings. For further details and proofs of these theorems, see \cites{MerkurjevSuslin82,OrlovVishikVoevodsky07,Voevodsky03a,Voevodsky03b,Weibel08,Weibel09}.

After the proof of the Bloch--Kato conjecture, the search for other properties of absolute Galois groups using Galois cohomology continued.  
Massey products are higher-degree analogs of the cup product that may be defined on the cohomology ring. The vanishing of the $n$-Massey products was conjectured by Min{\'a}{\v{c}}  and Tan in \cites{MT16,MT17b}. For the history of this conjecture, we refer the reader to \cites{HW19, MS23}. For some earlier developments leading up to this conjecture, see \cite{GLMS}. Related papers in Galois theory in a broader context are \cites{BCG23,HLMR24,PQ22}.

\subsection{Existence of decomposing fields}  \label{sec:finite-Bloch--Kato}
For each finite Galois extension $K/F$ and each element $\omega$ in $ H^m(\Gal(K/F), \fp)$, the existence of a decomposing field for $\omega$ --- a finite Galois extension $L/F$ that contains $K/F$ such that $\omega$ fully decomposes over $L$ --- follows immediately from the Bloch--Kato conjecture, or more specifically from the surjectivity of the map $\eta$ in (\ref{eq:norm.residue.homomorphism}). However, if we can establish the existence of this decomposing field directly by elementary methods using field theory and group cohomology, then such an approach will give a direct (elementary) proof of the surjectivity of $\eta$ without using the Bloch--Kato conjecture.  This follows from the following simple proposition. 
\begin{pro} \label{prop:finite-Bloch--Kato}
Let  $A_i^\bullet$,  $i \in (I , \prec)$ be a directed system of graded rings and let $A^\bullet := \underset{\longrightarrow}{\lim} \; A_i^\bullet$  be its direct limit. The following are equivalent for each $m \ge 1$.
\begin{enumerate}
    \item For each $i \in I$ and each $\omega \in A_i^m$, there exists $j \in I$ with $i \prec j$ such that $\omega$ fully decomposes in $A_j^\bullet$.
    \item  All elements of $A^m$ are fully decomposable.
\end{enumerate}
\end{pro}
The proof of this proposition is formal and straightforward, and is therefore omitted. The application of this proposition to our situation arises from the special case where 
\[A^\bullet = H^*(G_F, \fp) = \underset{\overset{\longrightarrow}{K}}{\lim} \; H^*(\Gal(K, F), \fp), \]
and $K$ runs over the directed system of all finite Galois extensions over $F$. Then condition (2) follows from the surjectivity of $\eta$ in  (\ref{eq:norm.residue.homomorphism}). The main context of the present paper is establishing (1) explicitly in many cases without invoking the Bloch--Kato conjecture.

\subsection{Kummer Theory} Although Kummer theory is hidden in the Bloch--Kato conjecture, it is worth recording separately as we use it frequently. Let $F$ be a field containing a primitive $p$-th root of unity for some fixed prime $p$, and let $G_{F}$ denote its absolute Galois group. Kummer theory gives a nice description of the group $H^{1}(G_{F}, \fp)$ where $G_{F}$ acts trivially on $\fp$. The starting point is the Kummer sequence
\[ 1 \lrar \mu_{p} \lrar \sep{F}^{\times} \lstk{x^{p}} \sep{F}^{\times} \lrar 1,\]
where the map $x^{p}$ raises elements to their $p$th powers, and the kernel of this map is the group of $p$-th roots of unity.  This is a short exact sequence of discrete $G_{F}$-modules which gives rise to a long exact sequence in 
Galois cohomology:
\[ 1 \rar \mu_{p}^{G_{F}} \lrar (\sep{F}^{\times})^{G_{F}} \lrar (\sep{F}^{\times})^{G_{F}} \lrar H^{1}(G_{F}, \mu_{p}) \lrar H^{1}(G_{F}, \sep{F}^{\times}) \lrar \cdots.\]
We identify $\mu_{p}$ with $\fp$ by choosing a particular root of unity $\xi_p \in \mu_p$ and then associating $\xi_p^x \in \mu_p$ with $x \in \fp$. Note that $G_{F}$ acts trivially on this, therefore $\mu_{p}^{G_{F}} \cong \fp$ as $G_F$-modules. 
Since $\sep{F}$ is Galois over $F$, we note that $(\sep{F}^\times)^{G_{F}} = F^{\times}$ and by Hilbert's theorem 90 the first cohomology group
$H^{1}(G_{F}, \sep{F}^{\times})$ is trivial. As a result, the first few terms of the above long exact sequence simplify to give
\[ 1 \rar \fp \lrar F^{\times} \lstk{x^{p}} F^{\times} \lrar H^{1}(G_{F}, \fp) \lrar 0,\]
which is equivalent to the Kummer isomorphism
\[ H^{1}(G_{F}, \fp) \cong F^{\times}/ (F^{\times})^{p}.\]
This isomorphism can be made explicit using the connecting homomorphism of the long exact sequence, which is as follows. A class  $[\alpha]$ in $ F^{\times}/ (F^{\times})^{p}$  corresponds to the cohomology class $\chi_\alpha \in H^{1}(G_{F}, \fp) \cong \Hom(G_{F}, \fp)$ which maps each $\sigma \in G_F$ to 
\[ \chi_\alpha(\sigma) = \frac{\sigma(\alpha^{1/p})}{\alpha^{1/p}}.\]
Observe that $\chi_\alpha$ does indeed have outputs in $\mathbb{F}_p$: since each $\sigma\in G_F$ permutes the roots of $x^p-\alpha=0$, the quotient of $\sigma(\alpha^{1/p})$ by $\alpha^{1/p}$ produces some root of unity, which we then associate with an element of $\mathbb{F}_p$ by using the identification $\mu_p \cong \fp$ from above. 

We note that  the degree-one part of the reduced Milnor $k$-theory is precisely the group $F^{\times}/ (F^{\times})^{p}$. Therefore, the Kummer isomorphism is exactly the degree one part of the Bloch--Kato conjecture. The Bloch--Kato map is obtained by taking the Kummer map and extending it to $H^{*}(G_{F}, \fp)$ using the universal property of the tensor algebra and realizing the Steinberg relations in the tensor algebra of $ F^{\times}/ (F^{\times})^{p}$ after Bass and Tate \cite{Bass-Tate}. 


The following result is well known among specialists (see, for example, the beginning of \cite{KLM} where it is implicitly stated and a proof is sketched). Therefore we omit its proof. 

\begin{pro}\label{pro:Cohology-F(p)}
    The inflation map 
    $$ H^m(\Gal(F(p)/F), \fp) \rightarrow H^m(G_F, \fp)$$ is an isomorphism 
    for any $m \geq 1$.
\end{pro}




\subsection{Characterizing classes in $H^2(G, \mathbb{F}_p)$}
Consider a profinite group $G$ that acts continuously on an abelian group $A$. That gives $A$ the structure of a $G$-module. The equivalence classes of extensions of $G$ by $A$ that induce the given $G$-module structure on $A$ are in 1-1 correspondence with elements of the cohomology group $H^2(G,A)$. In particular, when $A$ is the cyclic group of order $p$, which we identify with $\mathbb{F}_p$, and the action of $G$ on $\mathbb{F}_p$ is trivial,  classes in the group $H^2(G,\mathbb{F}_p)$ correspond to isomorphism classes of extensions of the form
\[1 \lrar \mathbb{F}_p \lrar H \lrar G \rar 1. \]
Moreover, under this bijection, the zero element in $H^2(G,\mathbb{F}_p)$ corresponds to the extension that splits.

Finally, we have the following $5$-term sequence in cohomology. If $N$ is a closed normal subgroup of a topological group $G$, then we have an exact sequence
\begin{equation}\label{eq:five.term.sequence}   0 
    \rar H^{1}(G/N, \fp) 
    \lstk{\text{inf}} H^{1}(G, \fp) 
    \lstk{\text{res}} H^{1}(N, \fp)^{G} 
    \lstk{\text{trg}} H^{2}(G/N, \fp) 
    \lstk{\text{inf}} H^{2}(G, \fp),
\end{equation}
where $\text{inf}$, $\text{res}$ and $\text{trg}$ are the inflation, restriction, and transgression maps, respectively, in 
cohomology.

\subsection{The Galois embedding problem}

We now move on to \emph{Galois embedding theory}, also used throughout. Suppose we are given a Galois extension $K/F$ with Galois group $G$ and a surjective homomorphism of groups
\begin{equation}\label{eq:pi} \pi \colon E \lrar G.\end{equation} 
When does there exist a Galois extension $L/F$ which contains $K/F$ such that $E \cong \Gal(L/F)$ and so that the restriction map 
\[ \mathrm{res}_{\Gal(L/F)}^{\Gal(K/F)} \colon \Gal(L/F) \lrar \Gal(K/F) \] 
coincides with the given map $\pi$? If such a field extension $L$ of $F$ exists, we say that the Galois embedding problem is solvable.

Here, we will discuss only a special case relevant to us. Namely, the case when the kernel of $\pi$ is a cyclic group of order $p$ (so can be identified with $\fp$). The embedding problem in question is called a Brauer-type embedding problem and has a particularly nice solution.  To explain this, consider the short exact sequence coming from the Galois embedding problem in (\ref{eq:pi}): 
\[1 \lrar \fp \lrar E \lstk{\pi} G \lrar 1.\]
We assume further that this extension is non-split and Frattini; that is, any subset $\mathcal{E}$ of $E$ with the property that $\pi(\mathcal{E})$ generates $G$ also generates $E$. The above short exact sequence corresponds to a cohomology class $\gamma$ in $H^2(G, \fp)$.

One can reformulate the embedding problem purely group-theoretically 
using the fact that we may consider $G$ as a quotient of the absolute Galois group $G_F$ of $F$. Any solution $L$ to the embedding problem in this setting comes with a surjection $G_F \rar \Gal(L/F) \cong E$ compatible with restrictions. In reverse, given a commutative diagram 
\[ \xymatrix{ & & & G_F \ar@{>>}[d] \ar@{>>}[ld] \\
1 \ar[r]& \fp \ar[r]& E \ar@{>>}[r] & G \ar[r] &1 }\]
the solution field $L$ is the fixed field of the kernel of the map $G_F \rar E$, and $\Gal(L/F)$ is isomorphic to $E$. 

For the purposes of this paper, it will be convenient to replace $G_F$ with other groups $\hat{G}$ that admit $G$ as a quotient. Here, we consider only groups $\hat{G} = \Gal(T/F)$ for some field extension $T$ of $K$; in doing so, we will limit ourselves to solutions to the Galois embedding problem that are given by intermediate fields between $F$ and $T$.  
Consider the inflation map 
$$\inf_G^{\hat G}: H^2(G, \fp) \rar H^2(\hat G, \fp).$$ The extension corresponding to $\inf(\gamma)$ is given by the top row of the diagram 
\begin{equation}\label{eq:inflation}
\xymatrix@=1em{
 1\ar[rr] & & \fp \ar[rr]  \ar@{=}[dd]  & & D \ar[rr] \ar[dd]  & & \hat G \ar[rr] \ar@{>>}[dd] & & 1    \\
   && && & \ar@{=>}[ul]& && \\
 1\ar[rr] & & \fp \ar[rr] & & E \ar[rr]^{\pi} & & G \ar[rr] & & 1    \\
},
\end{equation} 
where the right square is a pull-back square and $D$ is the fiber product of $E$ and $\hat G$ over $G$. We will use the theorem below throughout. 
\begin{thm}[{\cite[1.1]{hoechsmann}, see also \cite[Proposition 3.5.9]{NeukirchSchmidtWingberg}}] \label{thm:Hoechsmann}
    The embedding problem in (\ref{eq:pi}) has a solution that is an intermediate field between $F$ and $T$ if and only if $\inf_{\Gal(K/F)}^{\Gal(T/F)}(\gamma)=0$, or equivalently, the top row in (\ref{eq:inflation}) splits.
\end{thm}
The existence of the splitting of (\ref{eq:inflation}), together with the assumption that the bottom short exact sequence in the above diagram is Frattini,  tells us that $E$ is indeed a quotient of $\hat G = \Gal(T/F)$. Therefore, $E$ is isomorphic to $\Gal(L/F)$ for some Galois extension $L/F$ as desired in the above embedding problem. In short, $\inf_{\Gal(K/F)}^{\Gal(T/F)}(\gamma)$ in $H^{2}(\Gal(T/F), \fp)$ is the obstruction to the Galois embedding problem within the extension $T/F$.



We now turn towards some examples and conditions on the solvability of the embedding problem. 

\begin{Ex}\label{ex:embedding.problems}\ 
    \begin{enumerate} 
        \item Let $K/F$ be a cyclic extension of degree $p^n$ for some $n \in \mathbb{N}$ and let $\sigma$ be a generator of $G= \Gal(K/F) \cong C_{p^n}$. Consider the Galois embedding problem associated with $K/F$ and the extension 
        $$ 1 \rar \fp \rar C_{p^{n+1}} \rar C_{p^{n}} \rar 1.$$ 
        Then by \cite[2.4.4]{ledet} this embedding problem is solvable if and only if a primitive $p$-th root of unity is a norm in $K/F$.
        \item\label{ex:Heisenberg.solvability.via.norms} Let $p$ be an odd prime and consider the Heisenberg group $H$ that is non-abelian of order $p^3$ and exponent $p$. Then 
$$H = \left< u,v,w : u^p = v^p = w^p = 1, vu = uvw, wu = uw, wv = vw\right>.$$ Let $K = F\left( \sqrt[p]{a}, \sqrt[p]{b}\right)$, with $a,b \in F$, be a $C_p \times C_p$ Galois extension of $F$. The first factor $C_p$ arises by attaching $\sqrt[p]{a}$ and the latter factor by attaching $\sqrt[p]{b}$, respectively; i.e., we have elements $\sigma_a,\sigma_b \in \Gal(K/F)$ with $\sigma_a(\sqrt[p]{a})=\xi_p\sqrt[p]{a}, \sigma_a(\sqrt[p]{b})=\sqrt[p]{b}, \sigma_b(\sqrt[p]{a})=\sqrt[p]{a},$ and $\sigma_b(\sqrt[p]{b})=\xi_p\sqrt[p]{b}$ so that $\Gal(K/F)=\langle\sigma_a,\sigma_b:\sigma_a^p=\sigma_b^p=[\sigma_a,\sigma_b]=1\rangle$.  Consider the embedding problem given by $K/F$ and the extension 
$$1 \rar \fp \rar H \overset{\pi}\rar \Gal(K/F) \rar 1,$$
where the map $\pi$ is defined by $\pi(u)=\sigma_a$ and $\pi(v)=\sigma_b$, and $1 \in \fp$ is mapped to $w\in H$. This embedding problem is solvable if and only if $b$ is a norm in $F(\sqrt[p]{a})/F$ (see \cite[2.4.6]{ledet}). 

In fact, we can understand this embedding problem cohomologically. Using the section $s:\Gal(K/F) \to H$ given by $s(\sigma_a^i\sigma_b^j)=u^iv^j$, one can directly calculate that this group extension is equivalent to the one induced by the element $(b)\cup(a) \in H^2(\Gal(K/F),\fp)$. Indeed, one can compute $$s(\sigma_a^i\sigma_b^j)s(\sigma_a^k\sigma_b^\ell)=w^{jk}s(\sigma_a^{i+j}\sigma_b^{k+\ell}),$$ and a parallel computation will show that the action of $(b)\cup(a)$ on $(\sigma_a^i\sigma_b^j,\sigma_a^k\sigma_b^\ell)$ returns $w^{jk}$ as well.


For the case $p=2$, there is a similar embedding problem using the Dihedral group of order $8$ instead of the Heisenberg group. 
    \end{enumerate} 
\end{Ex} 

\section{Quotients of absolute Galois groups} \label{quotients}

In this section, we develop some results about quotients of the absolute Galois group. It should be noted that the main result of this section appeared in \cite{CEM}, where it was proved using abstract group cohomological techniques.  Here, we give a completely different proof using Galois embedding techniques. This not only keeps the paper self-contained but also gives new insight into the quotients. For ease of notation, we denote $\inf_{[n]}^{[m]} := \inf_{G_{F}^{[n]}}^{G_F^{[m]}}$ and $\inf_{[n]} := \inf_{G_{F}^{[n]}}^{G_F}$.

\begin{lem}\label{le:smallest.field.to.capture.first.cohomology}  The inflation map $\inf_{[2]} \colon H^1(G_F^{[2]}, \fp)  \lrar H^1(G_F, \fp)$ is a natural isomorphism.
\end{lem}

\begin{proof}
We remark that $H^1(G, \fp) \cong \Hom(G, \fp)$ for $G \in \{G_F, G_F^{[2]}$\}. The inflation map in question sends $f \colon G_F^{[2]} \rar \fp$ to $f \circ \pi \colon G_F \rar \fp$, where $\pi \colon G_F \twoheadrightarrow G_F^{[2]}$ is the quotient map. By construction of $F^{(2)}$, every element in $\Hom(G_F, \fp)$ lifts to some homomorphism in $\Hom(G_F^{[2]}, \fp)$ and so the inflation is surjective.  
The result now follows easily.
\end{proof}

\Cref{le:smallest.field.to.capture.first.cohomology} says that $F^{(2)}$ is the smallest field extension of $F$ that captures $H^1(G_F, \fp)$. This observation, although
elementary, is very crucial in the next theorem. To state the next theorem, we first need a definition.

\begin{defn}
Let $R_*$ be a connected graded $\fp$-algebra. A homogeneous element is said to be \emph{fully decomposable} if it can be generated by one-dimensional
classes, and is said to be \emph{not fully decomposable} otherwise.
 The \emph{decomposable part of $R_*$} is defined to be the $\fp$-subalgebra of $R_*$
that is generated by one-dimensional classes. This subalgebra will be denoted by $\Dec(R_*)$.
Similarly, a homogeneous element of $R_*$ is said to be \emph{reducible} if it can be generated by classes in lower degrees and is said to be \emph{irreducible} otherwise.
\end{defn}

We remark that the above aligns with \Cref{defn:decomposing-field} that characterizes decomposing fields of cohomology classes. 

\begin{thm}[{\cite[Theorem 6.6(a)]{CEM}}] \label{thm:G_F^{[3]}} 

Suppose that $L/F$ is Galois with $F^{(3)} \subseteq L$.
Then the restriction of the inflation map $\inf_{\Gal(L/F)}^{G_F}$ to $\Dec\left(H^*(\Gal(L/F),\fp)\right)$ is an isomorphism to $H^*(G_F, \fp)$. 

\end{thm}

\begin{proof}
The inflation map $\inf_{[2]} \colon H^*(G_F^{[2]}, \fp)  \lrar H^*(G_F, \fp)$ induced by the quotient map $G_F \rar G_F^{[2]}$ is shown in the above lemma to be an isomorphism in degree $1$.  By the Bloch--Kato conjecture, we know that $H^*(G_F, \fp)$ is generated by elements in degree $1$.
This tells us that the induced map $\inf_{[2]} \colon \Dec (H^*(G_F^{[2]}, \fp)) \rar H^*(G_F, \fp)$ is surjective.
The following commutative diagram of inflation maps
\[
\xymatrix{
\Dec\left( H^*(G_F^{[2]}, \fp)\right) \ar[r]^{\inf_{G_F^{[2]}}^{\Gal(L/F)}} \ar@{->>}[dr]_{\inf_{[2]}} & \Dec\left(H^*(\Gal(L/F), \fp)\right) \ar[d]^{\inf_{\Gal(L/F)}^{G_F}} \\
& H^*(G_F, \fp)
}
\]
then implies that the map $\inf_{\Gal(L/F)}^{G_F} \colon \Dec\left(H^*(\Gal(L/F), \fp)\right) \lrar H^*(G_F, \fp)$ is also surjective.

To see injectivity, we first claim that there is a map $\rho_L \colon k_*(F) \twoheadrightarrow  \Dec \; H^*(\Gal(L/F), \fp)$
that makes the following diagram commutative:
\[
\xymatrix{
& \Dec\left(H^*(\Gal(L/F), \fp)\right) \ar@{->>}[d]^{\inf_{\Gal(L/F)}^{G_F}} \\
k_*(F) \ar@{..>>}[ur]^{\rho_L} \ar[r]_{\cong\;\;\;} & H^*(G_F, \fp).
}
\]
Note that injectivity follows from the commutativity of the above diagram. So, it remains to establish the
existence of $\rho_L$. Furthermore, if we can show the existence of $\rho_3 :=\rho_{F^{(3)}}$, then we can create the desired $\rho_L$ as the composition of $\rho_3:k_*(F)\to\Dec\left(H^*(G_F^{[3]},\fp)\right)$ with $\inf_{[3]}^{\Gal(L/F)}:\Dec\left(H^*(G_F^{[3]},\fp)\right)\to\Dec\left(H^*(\Gal(L/F),\fp)\right)$. 

The existence of the (surjective) map $\rho_3$ will be guaranteed once we can show that the Steinberg relations $\{a, 1-a \} = 0$
also hold in $H^2(G_F^{[3]}, \fp)$. In other words, for any $a \ne 1$ in $F^\times$, we have to show that
the cup product of the classes $(a)$ and $(1 -a)$ is trivial in $H^2(G_F^{[3]}, \fp)$. In what follows, for an element $a \in F^\times$ we will write $[a]$ for the class it represents in $F^\times/F^{\times p}$, and will write $(a)$ for the element in the cohomology ring of $F$ corresponding to $a$. Of course, parentheses can also be employed in some other familiar ways, so the reader is encouraged to use context to infer whether or not an element is from a cohomology group.

Assume $p$ is an odd prime.
We then have to consider two cases.

\noindent
Case 1: $[a]$ and $[1 -a] $ are linearly dependent in $F^\times/F^{\times p}$.  Then we have the relation $[1 - a] = i[a]$ for some  $i$ in $\fp$,
which tells us that $(a) (1 - a) = (a)(a) i$. Since $p$ is odd, the square of any odd-degree class is zero in
$H^*(G_F^{[2]}, \fp)$.  In particular, $(a)(a) = 0$.

\noindent
Case 2: $[a]$ and $[1 -a] $ are linearly independent. 
By the assumption, the field $K= F\left( \sqrt[p]{a}, \sqrt[p]{1-a}\right)$ is Galois over $F$ with Galois group $C_p \times C_p$. We remark that 
$$  a = 1-(1-a) = \prod_{i=0}^{p-1} \left(1- \zeta_p^i \sqrt[p]{1-a} \right) = \mathrm{N}_{F\left( \sqrt[p]{1-a}\right)/F} \left(1-\sqrt[p]{1-a}\right). $$ This can be seen from factoring the polynomial $x^p-(1-a)$ into products of linear factors $(x-\zeta_p^i\root{p}\of{1-a})$ (for $i \in \{0,1,\cdots,p-1\}$) and plugging in $x=1$. (See also \cite[Section 8, Proof of Lemma 8.1]{Srinivas}.)


Since the previous displayed equation shows us that $a$ is a norm for the extension $F(\sqrt[p]{1-a})/F$, Example \ref{ex:embedding.problems}(\ref{ex:Heisenberg.solvability.via.norms}) tells us that the Galois embedding problem for $F(\sqrt[p]{a},\sqrt[p]{1-a})/F$ corresponding to the appearance of the Heisenberg group $H$ in the short exact sequence
\[1 \rar \fp \rar H  \stk{\pi} \Gal(K/F) \rar 1\] 
is solvable. (When $H$ is generated by $u,v,$ and $w$ as per Example \ref{ex:embedding.problems}(\ref{ex:Heisenberg.solvability.via.norms}), the map $\pi$ is defined by $\pi(u)=\sigma_{1-a}$ and $\pi(v)=\sigma_{a}$.)

We let $\hat L$ be a field which solves this embedding problem, so that $\Gal(\hat L/F) = H$ and the restriction $\Gal(\hat L/F) \to \Gal(K/F)$ is $\pi$. Also from Example \ref{ex:embedding.problems}(\ref{ex:Heisenberg.solvability.via.norms}) we have that the embedding problem over $\Gal(K/F)$ given by $\pi$ corresponds to $(a) \cup (1-a) \in H^2(\Gal(K/F),\fp)$, and by Theorem \ref{thm:Hoechsmann} we therefore have that the inflation of $(a) \cup (1-a)$ to $H^2(\Gal(\hat L/F),\fp)$ --- which we also denote by $(a) \cup (1-a)$ --- is trivial. Note in particular that we can express this element as a cup product of $1$-cycles in $H^2(\Gal(K/F),\fp)$, but that the product does not become trivial until we reach $H^2(\Gal(\hat L/F),\fp)$. Since $\hat L \subseteq F^{(3)}$, we then get $(a) \cup (1-a)$ vanishes in $H^2(G_F^{[3]},\fp)$, as desired.


The proof for $p=2$ is similar, but instead of the Heisenberg group, we have to use the dihedral group of order $8$ and a cyclic group of order $4$; see \cite[Prop.~2.3, 2.4]{MinacSpira96} for details. 
\end{proof}

\section{The second cohomology refinement of the Bloch--Kato conjecture}  \label{refinement}
\noindent

We begin by explaining the second cohomology refinement in a simple example. 
Let $F$ be a field such that $G_F(p) = \ints_p$, the ring of $p$-adic integers for some odd prime $p$. The existence of such a field can be seen as follows.
If $q$ is a power of some prime different from $p$ (so that $\fq$ admits
a primitive $p$-th root of unity), then the absolute Galois group of the field $\fq$ is $\widehat{\ints}$, the profinite completion of $\ints$. It has a Sylow-$p$ subgroup whose fixed field $F$ has the property $G_F(p) = \ints_p$. In fact, $F$ is the compositum of all cyclic extensions of degree coprime to $p$.

Now consider the filtration tower associated with the field $F$. For any $n$, the extension $F^{(n)}$ is cyclic over $F$ of degree $p^{n-1}$.  
We have
\[ G_F(p) \cong \ints_p \cong \invlim{i} \Gal(F^{(n)}/F) \cong \invlim{i} C_{p^{n-1}}. \]
By \cite[Thm.~7.3]{Carlson}, we know that $H^*(C_{p^{n-1}},\fp) \simeq \fp[X,Y_{n-1}]/(X^2)$, where $X$ is a degree $1$ element and $Y_{n-1}$ is a degree $2$ element that corresponds to a non-split extension of $C_{p^{n-1}}$ by $\fp$. Although there are $p-1$  distinct non-split group extensions of $C_{p^{n-1}}$ by $\mathbb{F}_p$, they all yield the same group: $C_{p^n}$. Specifically, suppose we identify generators $\sigma \in C_{p^{n-1}}$ and $\widetilde{\sigma} \in C_{p^n}$. A choice of a non-split extension of $C_{p^{n-1}}$ by $\mathbb{F}_p$ is equivalent to a choice of $1 \leq c \leq p-1$ with $\widetilde{\sigma}^p = \sigma^c$; for convenience we will choose $Y_{n-1}$ to correspond to the extension where $c=1$. For further details concerning the interpretation of multiples of $Y_{n-1}$ via the 
corresponding group extensions, see \cite[Ch.~4, Sec.~3, pp.~109--110]{MacLane}.

We now claim that $Y_{n-1}$ vanishes when inflated to $H^2(C_{p^n},\fp)$. Consider the non-split extension corresponding to $Y_{n-1}$ in $H^2(C_{p^{n-1}}, \fp)$: 
$$1 \rar \fp \rar C_{p^n} \rar C_{p^{n-1}} \rar 1. $$
The inflation of this class is given by the top horizontal extension in the diagram below, where the right commutative square is a pull-back square: 
\[
\xymatrix@=1em{
 1\ar[rr] & & \fp \ar[rr]  \ar@{=}[dd]  & & D \ar[rr] \ar[dd]  & & C_{p^n} \ar[rr] \ar[dd] & & 1    \\
   && && & \ar@{=>}[ul]& && \\
 1\ar[rr] & & \fp \ar[rr] & & C_{p^n} \ar[rr] & & C_{p^{n-1}}\ar[rr] & & 1    \\
}
\]
The top horizontal extension splits via the diagonal map $C_{p^n} \rar D$, so that $\inf(Y_{n-1}) = 0$. (Equivalently, there is a solution to the embedding problem associated to $Y_{n-1}$ in the $p$-cyclic field extension $F^{(n)}$ by \Cref{thm:Hoechsmann}.)

Hence if we consider the directed system for cohomology along the filtration tower, we have
\begin{align}
H^*(G_F, \fp) & \cong \dirlim{n} H^*(C_{p^{n-1}}, \fp) \nonumber  \\  
              & \cong \dirlim{n} \fp[X, Y_{n-1}]/(X^2) \ 
              \\
             & \cong \fp[X]/(X^2). \ \ \ \ \text{} \nonumber
\end{align}
Notice that the class $X$ survives because $X$ has degree one, and we know that
the cohomology map in degree one is an isomorphism. 

Using the terminology of \Cref{defn:decomposing-field}, the discussion in the above example can be summarized by saying that \emph{ ``$F^{(n)}$ is the minimal decomposing
field for the class $Y_{n-1}$.''}
In the following, we employ a similar idea as well as \Cref{thm:G_F^{[3]}} to prove \Cref{thm:main}. More precisely, we construct a minimal decomposing field for a general cohomology class in the filtration tower associated with $F$. 

Before we leave this example, we should remark that it is not necessary for a cohomology class to become trivial to fully decompose; see \Cref{sec:splitting.without.vanishing}.

\begin{proof}[Proof of \Cref{thm:main}]
 Let $\gamma$ be a class in $H^2(G_F^{[n]}, \fp)$ that is not fully decomposable. By the Bloch--Kato conjecture, we know $\inf_{[n]}(\gamma)\in H^2(G_F, \fp)$ is fully decomposable. Moreover, we have shown in \Cref{thm:G_F^{[3]}} that since $n \geq 3$, the fully decomposable class $\inf_{[n]}(\gamma)$ is the image of some fully decomposable class $\gamma_n$ in $H^2(G_F^{[n]}, \fp)$. By construction of $\gamma_n$, we have $\inf_{[n]}(\gamma-\gamma_n)$ vanishes in $H^2(G_F, \fp)$.
By \Cref{thm:Hoechsmann} the vanishing of $\inf_{[n]} (\gamma - \gamma_n) $  implies the solvability of the Galois embedding
problem given by the group extension
\[1 \rar \fp \rar Y \stk{\psi} G_F^{[n]} \rar 1 \]
corresponding to the class $\gamma -\gamma_n$ in $H^2(G_F^{[n]}, \fp)$.

Observe that $\gamma-\gamma_n$ is non-zero because $\gamma$ was assumed not fully decomposable. Therefore the group extension above is a nontrivial group extension of a pro-$p$ group by a cyclic extension of degree $p$. We claim that $\ker(\psi)$ is in the Frattini subgroup $\Phi(Y)$ of $Y$. To see this is true, recall that $\Phi(Y)$ is the intersection of all subgroups of $Y$ of index $p$. If $\ker(\psi)$ was not in $\Phi(Y)$, there would exist some subgroup $Z \subseteq Y$ of index $p$ such that $\ker(\psi) \not\subseteq Z$. But since $|\ker(\psi)|=p$, this would mean $\ker(\psi)\cap Z=\{1\}$. We also know that $\ker(\psi)$ is central in $Y$ (since $p$-groups act trivially on cyclic $p$-groups). But since $Z$ is a maximal subgroup, this would mean that $Y=Z \times \ker(\psi)$, in which case we would have $Z \simeq Y/\ker(\psi) \simeq G_F^{[n]}$. This, however, would mean the short exact sequence splits, contrary to the fact that $\gamma-\gamma_n\neq 0$.  

In this case it is easy to see that the embedding solution must be proper. Since $\ker(\psi) \subseteq \Phi(Y)$, any lift of generators of $G_F^{[n]}$ to $Y$ under $\psi$ will also generate $Y$. Indeed, \cite[Thm.~4.10]{Koch} tells us that if $\tilde G$ is a pro-$p$-group, then a collection generates $\tilde G$ if and only if the image of that collection in $\tilde G/\Phi(\tilde G)$ generates $\tilde G/\Phi(\tilde G)$. Because $\ker(\psi) \subseteq \Phi(Y)$ we see that $Y/\Phi(Y) \simeq G_F^{[n]}/\Phi(G_F^{[n]})$. As a consequence, if we lift any set of generators of $G_F^{[n]}$ to $Y$ under $\psi$, those lifts must also generate $Y$. Hence any homomorphism to $Y$ which solves our embedding problem must be surjective. This means that there exists a surjective homomorphism $G_F \to Y$ lifting the natural surjective homomorphism $G_F \to G_F^{[n]}$. By Galois theory, the fixed field of the kernel of this homomorphism onto $Y$ corresponds to some Galois extension $L/F$. (For a more complete treatment of the ideas in these paragraphs, the reader is encouraged to consult \cite[Statements 2.2 and 2.3, pg.~85]{hoechsmann}.)

In summary, there is a field extension $L$ of degree $p$ over $G_F^{[n]}$ that is  Galois over $F$, such that the restriction map $\Gal(L/F) \rar \Gal(G_F^{[n]}/F)$ is isomorphic to $\psi$. Another application of Theorem \ref{thm:Hoechsmann} --- this time using the fact that the above embedding problem is solvable by a field contained in the extension $L/F$ (namely $L$ itself) --- tells us that $\inf_{G_F^{[n]}}^{\Gal(L/F)}(\gamma - \gamma_n)$ is zero in 
$H^2(\Gal(L/F), \fp)$. 
Hence $\inf_{G_F^{[n]}}^{\Gal(L/F)}(\gamma)$ fully decomposes, as desired.
\end{proof}

We deduce the following corollary as a direct consequence of \Cref{thm:main}. A group-theoretic version of this result also appears in \cite[Proposition 5.6]{CEM}. 
\begin{cor}\label{cor:H2toH2dec} For any $n \geq 3$, every class $\gamma$ in $H^2(G_F^{[n]}, \fp)$ fully decomposes under the inflation map to $G_F^{[n+1]}$, i.e.,
$\inf_{[n]}^{[n+1]} (\gamma)$ belongs to $\Dec\left(H^2(G_F^{[n+1]}, \fp)\right)$.
\end{cor} 
In the following, we denote by $G_F^{(n)}$ the absolute Galois group $\Gal(F_{sep}/F^{(n)})$. Let $J_n = F^{(n)\times}/{(F^{(n)\times})}^p$ be the group of $p$-classes of $F^{(n)\times}$, and
let $\bar{J_n} = J_n^{G_F^{[n]}}$ be the submodule of $p$-classes invariant under the action of $G_F^{[n]}$. For $\gamma \in F^{(n)\times}$, we let $[\gamma]$ denote the $p$-class in $J_n$ represented by $\gamma$.
\begin{pro}\label{pro:Jn}
    Let $n \geq 3$. Then 
    \begin{enumerate}
        \item $\bar{J}_n = \ker\left[\inf_{[n]} \colon H^2(G_F^{[n]} , \fp) \lrar H^2(G_F, \fp)\right]$, and
        \item $F^{(n+1)} = F^{(n)} \left(\left(\bar{J_n}\right)^{1/p}\right)$.
    \end{enumerate}
\end{pro}

\begin{proof}
 (1)   The exact sequence
\[ 1 \rar G_F^{(n)} \rar G_F \rar G_F^{[n]} \rar 1\]
yields the fundamental 5-term sequence
\[ 0 \rar H^1(G_F^{[n]}, \fp) \stk{\inf} H^1(G_F, \fp) \stk{\text{res}} H^1(G_F^{(n)}, \fp)^{G_F^{[n]}} \stk{\text{tr}}  H^2(G_F^{[n]}, \fp) \stk{\inf} H^2(G_F, \fp).  \]
Since the first inflation map is an isomorphism, the kernel of the second inflation map is
\[H^1(G_F^{(n)}, \fp)^{G_F^{[n]}},\]
which, by Kummer theory, is isomorphic to
\[ \left(F^{(n)\times}/ {(F^{(n)\times})}^p \right)^{G_F^{[n]}} = \bar{J}_n. \]

(2)
By our construction, for each $n$ we have the following short exact sequence:
\[ 1 \lrar G_{F}^{(n)}/G_{F}^{(n+1)} \lrar G_{F}^{[n+1]} \lrar G_{F}^{[n]} \lrar 1. \] By equivariant Kummer theory (see, for example, \cite[Sec.~1]{WaterhouseEquivariantKummer}\footnote{Waterhouse discusses equivariant Kummer theory in the case of finite Galois extensions, but his proof is valid when $V$ is an infinite dimensional module and when we use the infinite Kummer pairing and Pontryagin duality between a compact group and its discrete dual (see \cite[Ch.~6, Sec.~2]{ArtinTate}).}), there is a trivial  $\fp[G_F^{[n]}]$-module $V \subseteq F^{(n) \times}/(F^{(n)\times})^p$ (i.e., an $\fp$-subspace of $\bar J_n$) such that the field $F^{(n+1)}$ is $F^{(n)}(k^{1/p}:k \in V)$. Indeed, we must have $V = \bar{J_{n}}$ because $G_{F}^{[n+1]}$ is the largest central extension of $G_{F}^{[n]}$ by an elementary $p$-abelian subgroup.
Therefore $F^{(n+1)} = F^{(n)}\left((\bar{J_{n}})^{1/p}\right)$, as desired.
\end{proof}


\begin{thm}\label{thm:F^(n+1)minimal-decomposing}
    Let $n \geq 3$. The field $F^{(n+1)}$ is a field minimal with respect to inclusion that is a decomposing field for all $\gamma \in H^2(G_F^{[n]}, \fp)$. More precisely, there is no field $L$ so that 
    \begin{enumerate}
        \item $F^{(n)} \subset L \subsetneq F^{(n+1)}$, 
        \item $L$ is Galois over $F$, and 
        \item all $\gamma \in H^2(G_F^{[n]},\fp)$ fully decompose over $L$. 
    \end{enumerate}
\end{thm}
\begin{proof}
    Let $L$ be a field extension of $F^{(n)}$ that is Galois over $F$ so that any cohomology class $\gamma \in H^2(G_F^{[n]},\fp)$ fully decomposes over $L$. Let $K/F^{(n)}$ be any degree $p$ extension of $F^{(n)}$ that is Galois over $F$. We know that $K$ is defined by some nonzero element $\overline{\theta} \in \bar{J}_n$ for $\theta \in F^{(n)\times}$. This uniquely determines an extension 
    \[ 1 \rightarrow \fp \rar  \Gal(K/F) \rightarrow G_F^{[n]} \rightarrow 1, \]
    and we know this extension must be non-split since otherwise we would have $K \subseteq F^{(n)}$. Hence the extension
     defines a non-trivial class $\gamma \in H^2(G_F^{[n]},\fp)$. Theorem \ref{thm:Hoechsmann} tells us that $\inf_{[n]}(\gamma)=0$ since the embedding problem corresponding to $\gamma$ is solvable in $G_F$.


    It remains to show that $K \subseteq L$. From our hypothesis (3), we see that $\inf_{G_F^{[n]}}^{\Gal(L/F)}(\gamma) \in \Dec\left(H^*(\Gal(L/F),\fp)\right)$. But since $\inf_{\Gal(L/F)}^{G_F}(\inf_{G_F^{[n]}}^{\Gal(L/F)}(\gamma)) = \inf_{[n]}(\gamma) = 0$, Theorem \ref{thm:G_F^{[3]}} implies that $\inf_{G_F^{[n]}}^{\Gal(L/F)}(\gamma)=0$. We consider the Brauer embedding problem
    \[ \xymatrix{ & & & \Gal(L/F)\ar[d] \ar@{.>}[ld]_\omega\\
    1 \ar[r] & \fp \ar[r] &\Gal(K/F) \ar[r] & G_F^{[n]} \ar[r] & 1}
    \]
    where the bottom extension is given by $\gamma$. By \Cref{thm:Hoechsmann}, there is some intermediate field $F^{(n)} \subseteq K' \subseteq L$ that is Galois over $F$ so that there is an isomorphism 
    \[ \pi: \Gal(K'/F) \rightarrow \Gal(K/F)\]
    making the following diagram commute
    \[ \xymatrix{ 1 \ar[r] & \fp \ar@{=}[d]\ar[r] & \Gal(K'/F) \ar[d]^\pi \ar[r] &G_F^{[n]}\ar@{=}[d] \ar[r] & 1 \\ 1 \ar[r] & \fp \ar[r] & \Gal(K/F) \ar[r] & G_F^{[n]} \ar[r] & 1 }  \]
    and $K'/F$ solves the embedding problem above. We recall that $K= F^{(n)}\left( \sqrt[p]{\theta}\right)$ and thus by \cite[Lemma A.1.1]{JLY2002}, there exists $f \in F^\times$ so that the field $K' = F^{(n)}\left( \sqrt[p]{f\theta} \right) $. By assumption $F^{(n)} \supset F^{(2)}$ and $\sqrt[p]{f} \in F^{(2)}$ so that $K' = K$ and therefore $K$ is contained in $L$. \\
    Because $K$ is any cyclic $p$-extension of $F^{(n)}$ of degree $p$ and $F^{(n+1)}$ is the compositum of all such extensions, we see that $F^{(n+1)} \subseteq L$. 
\end{proof}

The above discussion leads to another interpretation of the filtration tower associated to $F$: $F^{(n+1)}$ is the field that simultaneously fully decomposes all degree two cohomology classes on $F^{(n)}$. As a result, $F(p) = \bigcup_{n\geq 1} F^{(n)}$ is the smallest field extension of $F$ that contains $F^{(3)}$ that is ``closed under decomposition of cohomology classes" (see also \cite[Theorem 8.5]{CEM}). 

\section{An example of splitting without vanishing}\label{sec:splitting.without.vanishing}

After the above results, it is natural to be interested in examples of elements which fully decompose without becoming trivial. In fact, there are many such examples. Here we provide an interesting example of an element that is not fully decomposable in the cohomology algebra of a certain finite Galois extension $K$ over $F$ with $\FF_2$-coefficients in degree four which fully decomposes without vanishing in a larger extension $L/F$.  Indeed this example can be viewed as the implementation of a widely applicable strategy for producing elements that fully decompose without vanishing under inflation: employing the known vanishing of a particular cohomology element under inflation to manufacture non-trivial decompositions.

\begin{Ex}\label{example}
    Suppose that $\text{char}(F) \neq 2$, and  we have a tower of Galois extensions $F \subset K \subset N \subset L$
    where:
    \begin{enumerate} 
    \item $\Gal(K/F) \cong Q_8 \times \left( \ZZ/4\right)^4$, where 
\[ Q_8=\langle g,h\mid g^4=1,\ g^2=h^2,\ gh=h^3g\rangle, \]
    is the quaternion group of order $8$; 
    \item $\Gal(N/F) \cong \left( \ZZ/4 \rtimes \ZZ/4\right) \times \left( \ZZ/4\right)^4$ with generators $\gamma$ and $\delta$ of the first and second factor, respectively, satisfying $\gamma \delta = \delta^{-1}\gamma$;  
    \item $\Gal(L/F) \cong \left( \ZZ/8 \rtimes \ZZ/8\right) \times \left( \ZZ/4\right)^4$ with generators $a$ and $b$ of the first and second factor, respectively, satisfying $a b = b^3 a$;
    \item the restriction of Galois actions from $N$ to $K$ induces the homomorphism 
    $$\varphi: \ZZ/4 \rtimes \ZZ/4\longrightarrow Q_8$$
    with kernel $\{1, \gamma^2\delta^2\}$ so that $\phi(\gamma) = i$ and $\phi(\delta) = j$ (for more details, see e.g. \cite{Conrad} or \cite{GSS95}); and
   \item the restriction from $L$ to $N$ gives a natural surjective homomorphism 
   $$\theta: \ZZ/8 \rtimes \ZZ/8 \longrightarrow \ZZ/4 \rtimes \ZZ/4 $$
   with $\theta(a) = \gamma$ and $\theta(b) = \delta$. 
   \end{enumerate}
   We first show that such a field extension exists. A well-known argument of Artin (see \cite[Ch.~2, Sec.~G]{Artin}) tells us that for any field $k$ and any finite group $G$, there exists some extension $K/F$ with $k \subseteq F$ and so that $\Gal(K/F) \simeq G$. 
   In particular, we can find $L/F$ satisfying (3) above. Using the homomorphisms $\theta$ and $\varphi$ we can define a tower of field extensions such that (1)--(5) hold.

From \cite[Prop.~3.5.5]{B} we have that $H^*(\mathbb{Z}/4,\FF_2)$ is generated by $X$ and $Y$ with $\deg(X)=1$ and $\deg(Y)=2$ and subject to the single relation $X^2=0$. Hence by K\"unneth's theorem \cite[Page 67, Remark]{B}, the cohomology algebra decomposes as 
   \begin{align*}
       H^*\left( \Gal(L/F), \FF_2\right) &\cong H^*\left( \left( \ZZ/8 \rtimes \ZZ/8 \right) \times \left( \ZZ/4\right)^4, \FF_2\right)\\
       &\cong H^*\left( \left( \ZZ/8 \rtimes \ZZ/8\right), \FF_2\right) \otimes_{\FF_2} \left( \bigwedge (X_1, X_2,X_3,X_4) \otimes_{\FF_2} \FF_2[Y_1, Y_2, Y_3, Y_4] \right),
   \end{align*}
    where $\bigwedge(X_1, X_2,X_3,X_4)$ is the exterior algebra of $\FF_2$ with $\deg(X_i) =1$ and $\deg(Y_i) =2$ for $1 \leq i \leq 4$. We obtain similar tensor decompositions for $H^*\left( \Gal(N/F), \FF_2\right)$ and $H^*\left( \Gal(K/F),\FF_2\right)$. The cohomology of $Q_8$  (see \cite[Ch.~4, Lem.~2.19, pg.~129]{AM04}) is 
    $$H^*\left( Q_8, \FF_2\right) = \FF_2\left[ X,Y,Z\right] / \left( X^2 + Y^2 + XY, X^2Y + XY^2\right),$$ where $\deg(X) = \deg(Y)=1$ and $\deg(Z)= 4$.
    We observe that $X^2+Y^2 = XY$ implies that $X^3 + XY^2 = X^2Y$ and thus by the second relation $X^3 =0$. Similarly, we see that $Y^3 = 0$. Therefore, we can take the following generators
    \begin{center}
    \begin{tabular}{c|c}
        degree &  generators \\ \hline 
        $0 + 4k$  & $Z^k1$ \\
        $1 + 4k$ & $Z^kX,Z^kY$\\
        $2 + 4k$ & $Z^kX^2,Z^kY^2$\\
        $3 + 4k$ & $Z^kX^2Y$
    \end{tabular}\end{center}
    for any $k \in \ZZ_{\geq 0}$. Thus the cohomology of $H^*\left(\Gal(K/F),\FF_2\right)$ contains the following element that is not fully decomposable:
    $$W = Z + X_1X_2X_3X_4.$$
    By the computations in Proposition A.1 (in the appendix), we have $\inf_{Q_8}^{\ZZ/8 \rtimes \ZZ/8} (Z) = 0$ and $\inf_{Q_8}^{\ZZ/4 \rtimes \ZZ/4} (Z) \neq 0$. Therefore
    $$ \inf_{\Gal(K/F)}^{\Gal(L/F)}(W) = X_1 X_2 X_3 X_4 $$
   is a nontrivial element that fully decomposes, but does not vanish, over $L$.
\end{Ex}

\begin{rem}
An additional field-theoretic example can be constructed in a way similar to Example \ref{example} using Proposition A.1 from the appendix. This example will yield a smaller Galois extension with the Galois group isomorphic to $H \times \mathbb{Z}/2$ (in the notation of Proposition A.1). Thus this Galois group will have order $64$.
\end{rem}

\section{Refining Bloch--Kato: Superpythagorean and $p$-rigid fields} \label{computation}

In this section we  provide a proof of \Cref{thm:example}. We do this by studying the decomposition of cohomology classes of degree $2$ associated to the absolute Galois groups of so-called superpythagorean and $p$-rigid fields. We will soon see that these fields have associated Galois groups of the form $\mathbb{Z}_2^I \rtimes C_2$ (in the former case) or $\mathbb{Z}_p^I \rtimes \mathbb{Z}_p$ (in the latter case), where $p$ is prime and $I$ is some indexing set.  Since these groups are profinite, their cohomology is determined by the cohomology of their finite quotients; the bulk of this section focuses on computing the cohomology of these finite quotients.

Our computation of these cohomology rings will illustrate the refinement of the Bloch--Kato conjecture proven earlier by observing directly how classes that are not fully decomposable become fully decomposable under inflation. It is also worth noting that the results of this section are developed independently of the Bloch--Kato conjecture.

We begin by defining the relevant groups, starting with the case when $p=2$. For $n,d \in \mathbb{N}$, define $A_2(n;d)$ as the semidirect product $A_2(n;d)=(C_{2^n})^d\rtimes C_2,$ where the action in the semidirect product is given by $xyx^{-1}=y^{-1}$ for $y \in (C_{2^n})^d, x \in C_2$. For an index set $I$, we will then be interested in the group $$\mathbb{Z}_2^I \rtimes C_2 = 
\varprojlim_{\substack{
  J \subseteq I \\
  |J|<\infty
}}
\ilim{n}{} A_2(n;|J|),
$$
where the interior inverse limit is taken with respect to the maps $A_2(n+1;|J|)\to A_2(n;|J|)$  induced by the natural surjection $C_{2^{n+1}} \to  C_{2^{n}}$ on the first coordinate and the identity map on the second.

When $p>2$ we have more actions to consider when constructing our semidirect product, so we will need an additional parameter in our definition. For $n,d \in \mathbb{N}$ and $k \in \mathbb{N}\cup \{\infty\}$, define $B_p(n;d,k) = (C_{p^n})^d\rtimes C_{p^n},$
where this time the action is given by $xyx^{-1}=y^{p^k+1}$ for all $y \in (C_{p^n})^d$ and $x$ a fixed generator of $C_{p^n}$ (the outer factor in our semidirect product). By convention, when $k=\infty$ the action is given by  $xyx^{-1}=y$, in which case the semidirect product is indeed the direct product. For an index set $I$ and a choice of $k \in \mathbb{N}\cup\{\infty\}$, we will then be interested in the group $$\mathbb{Z}_p^I \rtimes \mathbb{Z}_p = 
\varprojlim_{\substack{
  J \subseteq I \\
  |J|<\infty
}}
\ilim{n}{} B_p(n;|J|,k),
$$
with the inverse limits defined as before.

We now define the fields of interest. A field $F$ is said to be Pythagorean if every sum of two squares (hence any number of squares) in $F$ is a square. Observe that $\mathbb{R}$ and $\mathbb{C}$ are Pythagorean, but $\mathbb{Q}$ is not. Fields of characteristic $2$ are Pythagorean because in those fields we have $(x+y)^2 = x^2+y^2$.  We say $F$ is formally real if $-1$ is not a sum of squares in $F$. These are precisely the fields that admit an ordering, making them into an ordered field. $\mathbb{R}$ and $\mathbb{Q}$ are formally real. A more interesting example is the field $\mathbb{R}(x)$ of rational functions. If $r(x) \in \mathbb{R}(x)$ takes the form $$r(x) = \frac{a_nx^n+\cdots+ a_0}{b_nx^n + \cdots + b_0}$$ then we say that $r(x)$ is positive if  $a_nb_n > 0$ in $\mathbb{R}$. Further examples of orderings can be found in \cite[pgs.~4,82--85]{KS}. A superpythagorean field is a formally real field in which every subgroup of index $2$ in $F^{\times}$ not containing $-1$ defines an ordering on $F$.   

A superpythagorean field is necessarily Pythagorean (see the footnote in \cite[pg.~45]{Lam}). An example of a superpythagorean field is the formal power series ring over $\mathbb{R}$; this field is typically denoted $\mathbb{R}((x))$. More generally, the field of iterated power series $\mathbb{R}((x_1))((x_2))\cdots((x_n))$ is superpythagorean.


When $F$ is superpythagorean, the Galois group $G_F(2)$ of the maximal $2$-extension of $F$ is isomorphic to $\mathbb{Z}_2^I \rtimes C_2$. (See, for example, \cite[Cor.~4.10]{MRT16} in the case where $G_F(2)$ is finitely generated.)

Fields that have a similar Galois structure are $p$-rigid fields, which are an important class of fields in Galois theory. 
For the purposes of this paper, when we discuss $p$-rigid fields we will focus only on the case where $p$ is an odd prime. Readers interested in learning more about $2$-rigid fields can consult \cite{BCW}. 

Let  $F^p$ denote the collection of elements in $F$ that are $p$th powers. An element $a$ in $ F\smallsetminus F^p$ is said to be $p$-rigid if the image of the norm map $F(\sqrt[p]{a})\rightarrow F$
is contained in $\bigcup_{k=0}^{p-1}a^kF^p$.
We say that $F$ is $p$-{\it rigid} if all of the elements of $F\smallsetminus F^p$ are $p$-rigid. It can be shown that the field $\mathbb{Q}_l(\zeta)$ is $p$-rigid whenever $\zeta$ is a primitive $p$-th root of unity and $l$ is a prime different from $p$. 

For our purposes, $p$-rigid fields are of particular interest because the Galois group $G_F(p)$ of the maximal $p$-extension of a $p$-rigid field (when $p$ is odd and $F$ contains a primitive $p$th root of unity) is a group of the form $\mathbb{Z}_p^I \rtimes \mathbb{Z}_p$, as defined above. (See \cite[Thm.~4.10]{CMQ} in the case where $p$ is odd.)


Thus we see that when $F$ is superpythagorean or $p$-rigid, the cohomology groups are given as direct limits along the inflation maps:
$$
H^*(G_F(p),\FF_p)=\left\lbrace\begin{array}{ll}
\dlim{} H^*(A_2(n;d),\FF_2) & \text{if $F$ is superpythagorean and $p=2$}\\
\dlim{} H^*(B_p(n;d,k),\FF_p) & \text{if $F$ is $p$-rigid and $p$ is odd}.
\end{array}\right.
$$  

Since the cohomology of $G_F(p)$ is isomorphic to the cohomology of $G_F$ (by Proposition \ref{pro:Cohology-F(p)}), this allows us to approach the proof of Theorem \ref{thm:example} by first studying the cohomology rings $H^*(A_2(n;d),\FF_2)$ and $H^*(B_p(n;d,k),\FF_p)$.

\subsection{Cohomology of metacyclic groups}\label{subsec:cohomology.of.metacyclic.groups}
Ultimately we aim to prove Theorem \ref{thm:example}, but to do this we first focus on the case where the parameter $d$ defining the relevant Galois groups is $1$. Before diving into this proof, however, we make some observations.

Under the assumption that $d=1$, the groups $A_2(n;1)$ and $B_p(n;1,k)$ are metacyclic, i.e., extensions of cyclic groups. Fortunately, the mod-$p$ cohomology ring of metacyclic groups is completely determined by Huebschmann in \cite{H}. We start by describing the cohomology ring of certain types of metacyclic groups.  When discussing the cohomology of groups, it often will be convenient to identify a cyclic group $C_q$ with the additive group $\ZZ/q$.

We restrict our attention to split metacyclic groups, which sit in an extension of the form
\begin{equation}\label{exP}
1\rightarrow \ZZ/{p^n}\rightarrow P\rightarrow \ZZ/{p^m}\rightarrow 1
\end{equation}
where $n\geq m>0$. The generators of the quotient and the kernel groups will be denoted by $x$ and $y$, respectively. The action of $x$ on $y$ is given by $xyx^{-1}=y^t$ such that $t^{p^m}\equiv 1 \mod p^n$ and $t\equiv 1 \mod p$. In general, a 
metacyclic group $P$ 
admits
a presentation of the form
$$
P = \Span{x,y|\; y^r=1,\;x^s=y^f,\; xyx^{-1} = y^t},
$$
where $r>1$, $s>1$, $t^s \equiv 1 \mod r$, and $tf\equiv f \mod r$. The case of interest here corresponds to the split situation with
$f=0$, $r=p^n$, and $s=p^m$. 

We recall that the cohomology ring of a cyclic $p$-group is 
$$
H^*(\ZZ/{p^{m}} , \fp)= \begin{cases} 
\bigwedge(\omega_{x})\otimes \FF_p[c_{x}] & \text{if } p \text{ is odd or } p=2 \text{ and } m>1 \\
 \FF_2[\omega_x] & \text{if }  p=2 \text{ and } m=1
\end{cases}
$$  
where  $\omega_x$ and $c_x$ are in degree-one and degree two, respectively. Under the identification  $H^1(\ZZ/{p^m}, \fp)\cong \Hom(\ZZ/{p^m},\ZZ/p)$, the class $\omega_x$ corresponds to the mod $p$ reduction map $\ZZ/{p^m}\rightarrow \ZZ/p$ defined by $1\mapsto 1$. To describe the class $c_x$, consider the exact sequence of abelian groups $1\rightarrow \ZZ/p\rightarrow \ZZ/{p^{m+1}}\rightarrow \ZZ/{p^m}\rightarrow 1$. This induces a long exact sequence in cohomology 
$$
 \cdots H^1(\ZZ/p^m,\ZZ/p^{m+1})\rightarrow  H^1(\ZZ/p^m,\ZZ/p^m) \stackrel{\beta_m}{\rightarrow} H^2(\ZZ/p^m,\ZZ/p) \cdots
$$
where the map $\beta_m$ is called the $m$-th Bockstein. Now the class $c_x$ is the image of the identity map in $H^1(\ZZ/p^m,\ZZ/p^m)\cong \Hom(\ZZ/p^m,\ZZ/p^m)$ under the Bockstein homomorphism. See \cite[Section 6.2]{M} for further details. 

The cohomology ring of the metacyclic group $P$ is calculated in \cite[Theorem B]{H}. This theorem covers metacyclic groups where all the numbers $r$, $s$, $t-1$, $(t^s-1)/r$, $f$, $(t-1)f/r$ are divisible by $p$. Our group $P$ satisfies these conditions. As a notational convenience we do not distinguish the  classes $\omega_x$ and $c_x$ from their image  under the inflation map $\inf_{\ZZ/p^m}^P:H^*(\ZZ/p^m,\FF_p)\rightarrow H^*(P,\FF_p).$

\begin{thm}[{\cite[Thm.~B]{H}}]\label{thm:Huebschmann}The cohomology spectral sequence of the extension describing $P$ collapses on the $E_2$-page, and as an $H^*(\ZZ/p^m,\fp)$-module we have
$$ H^*(P, \fp)\cong H^{*}(\ZZ/p^m, \fp)\otimes H^{*}(\ZZ/p^n, \fp).$$
Moreover, there are classes $\omega$ and $c$ of degree-one and two, which restrict to the classes $\omega_y$ and $c_y$, respectively.
\begin{enumerate}[(i)]
\item For $p$ odd, as a graded commutative algebra we have
$$
H^*(P, \fp)=\bigwedge(\omega_x,\omega)\otimes \FF_p[c_x,c].
$$
\item For $p=2$, the cohomology ring is generated by $\omega_x$, $c_x$, $\omega$, $c$  as a graded commutative algebra  subject to the relations
\begin{align*}
\omega_x^2&=2^{m-1}c_x\\
\omega^2&=2^{n-1}(c_x+c)+\frac{t-1}{2}\omega_x\omega.
\end{align*}
\end{enumerate}
\end{thm}

\begin{rem}
The result \cite[Thm.~E(2)]{H} includes an error that was corrected in \cite{HuebEr}.  The above result comes from Theorem B of \cite{H} which was not affected by this error.
\end{rem}

Assume $n\geq 2$. 
For $A_2(n;1)$ we have $t=-1$ and $m=1$, and therefore
$$
H^*(A_2(n;1), \FF_2)= \FF_2[\omega,c,\omega_x]/(\omega^2+\omega_x\omega),
$$
where we have used the relation $c_x=\omega_x^2$. Since $n \geq 2$, the summand $2^{n-1}(c_x + c)$ vanishes.

Next consider $B_p(n;1,k)$; recall this means $p$ is an odd prime. We have $n=m$ and $t=p^k+1$ for some $k>0$. 
We have
$$
H^*(B_p(n;1,k), \fp)=\bigwedge(\omega_x,\omega)\otimes \FF_p[c_x,c].
$$

Let $P(n)$ stand either for $A_2(n;1)$ or $B_p(n;1,k)$. We want to study the inflation map
$$
\inf_{P(n)}^{P(n+1)}:H^*(P(n),\fp) \rightarrow H^*(P(n+1), \fp)
$$ induced by the natural projection $P(n+1) \rightarrow P(n)$. To distinguish the classes belonging to $P(n)$  for different $n$'s, we will use the notation $\theta(n)$ when we want to indicate that a cohomology class $\theta$ belongs to the $n$-th group. The one-dimensional classes $\omega_x(n)$ and $\omega(n)$ map to $\omega_x(n+1)$ and $\omega(n+1)$ under the inflation map, respectively. 

\Pro{\label{dec}Under the inflation map $\inf_{P(n)}^{P(n+1)} \colon H^*(P(n), \fp)\rightarrow H^*(P(n+1), \fp)$, two-dimensional classes fully decompose into a product of one-dimensional classes.}  
\begin{proof}
We begin with the case $P(n) = B_p(n) := B_p(n;1,k) = (\ZZ/p^n)\rtimes \ZZ/p^n$. Consider the exact sequence 
\[ \xymatrix{ 1 \ar[r] &\ZZ/p^{n} \ar[r] & B_p(n) \ar[r]& \ZZ/p^n \ar[r]& 1.} \]
Recall that the cohomology ring $H^*(B_p(n),\fp)$ has four generators. The one-dimensional class $\omega_x(n)$ and the two-dimensional class $c_x(n)$ are inflated from the quotient group. (The corresponding classes in the cohomology ring of the quotient group and their inflated images are denoted by the same symbols.) The other two generators, the one-dimensional class $\omega(n)$ and the two-dimensional class $c(n)$, restrict to the classes $\omega_y$ and $c_y$ of the kernel group, respectively.

A general principle we will use when computing the images of inflation maps is the following.
As observed in Section~\ref{refinement}, under the inflation map
\[
\inf_{\ZZ/p^n}^{\ZZ/p^{n+1}} \colon H^*(\ZZ/p^n,\fp)\longrightarrow H^*(\ZZ/p^{n+1},\fp),
\]
the image of the extension class is zero. To compute the inflation map associated to the middle horizontal arrow in diagram (\ref{eq:B_p diagram}), we will apply this observation to the inflation maps induced by the maps between the kernels and the quotients:
\begin{equation}\label{eq:B_p diagram}
\vcenter{
    \xymatrix{
\ZZ/p^{n+1}  \ar@{->}[d] \ar@{->}[r] & \ZZ/p^n \ar@{->}[d] \\
B_p(n+1)\ar@{->}[d] \ar@{->}[r] & B_p(n)\ar@{->}[d] \\
\ZZ/p^{n+1}  \ar@{->}[r] & \ZZ/p^n.  \\
}
}
\end{equation}
We then exploit the commutativity of the squares in cohomology induced by the diagram above.

For $P(n) = B_p(n)$, the classes $c_x(n)$ and $c(n)$ do not fully decompose, however, we will show that after applying inflation, they fully decompose. First, $\inf_{B_p(n+1)}^{B_p(n)}c_x(n)=0$ because the class $c_x(n)$ can be identified with the extension class of $\ZZ/p^{n+1}$, and hence maps to zero under inflation. 
It remains to determine the inflation of $c(n)$. 
We show that the class $\inf_{B_p(n)}^{B_p(n+1)}c(n)$ fully decomposes into a product of one-dimensional classes. 
In degree two $H^*(B_p(n+1), \fp)$ is generated by $c(n+1)$, $\omega_x(n+1)\omega(n+1)$, and $c_x(n+1)$. Let us  write
$$\inf_{B_p(n)}^{B_p(n+1)}c(n)=\alpha_1 c(n+1)+\alpha_2 \omega_x(n+1)\omega(n+1)+\alpha_3c_x(n+1)$$
for some $\alpha_i$ in $\FF_p$. Our goal is to show that $\alpha_1=\alpha_3=0$.
Let $\res_{B_p(n)}^{\ZZ/p^n}:H^*(B_p(n), \fp)\rightarrow H^*(\ZZ/p^n, \fp)$ denote the restriction map induced by the inclusion $\ZZ/p^n\rightarrow B_p(n)$. 
Consider the following commutative diagram induced by the top square in diagram (\ref{eq:B_p diagram}):
\begin{equation}\label{eq:inf res square}
\vcenter{
\xymatrix@=1em{
H^*(B_p(n),\fp) \ar[rr]^{\inf_{B_p(n)}^{B_p(n+1)}} \ar[dd]_{\res_{B_p(n)}^{\ZZ/p^n}} && H^*(B_p(n+1),\fp) \ar[dd]^{\res_{B_p(n+1)}^{\ZZ/p^{n+1}}}\\ \\
H^*(\ZZ/p^n, \fp) \ar[rr]^{\inf_{\ZZ/p^n}^{\ZZ/p^{n+1}}} && H^* (\ZZ/p^{n+1}, \fp).
}
}
\end{equation}
We will use the fact that $\omega_x(n)$ and $c_x(n)$ restrict to $0$ since they are inflated from the quotient. 
Note that the composite $\res_{B_p(n)}^{\ZZ/p^n} \circ \inf_{\ZZ/p^{n}}^{B_p(n)}$
vanishes in all positive degrees, since it is induced by the composition of the kernel inclusion with the quotient homomorphism, which is trivial.
Moreover, $\res_{B_p(n)}^{\ZZ/p^n} (c(n))$ coincides with $c_y(n)$, which inflates to zero since it is the extension class.
Using the commutativity of  
this
diagram for $c(n)$ together with these observations, we calculate that 
$$\alpha_1 c_y(n+1) = \res_{B_p(n+1)}^{\ZZ/p^{n+1}} \inf_{B_p(n)}^{B_p(n+1)} ( c(n)) =  \inf_{\ZZ/p^n}^{\ZZ/p^{n+1}}\res_{B_p(n)}^{\ZZ/p^n} (c(n)) = 0, $$
and so $\alpha_1 =0$.  
Next, we show that $\alpha_3=0$. Consider the splitting $s_n$ of the surjection $B_p(n)\rightarrow \ZZ/{p^n}$ defined by $s_n(1)=(0,1)$. The map $s_n^*$ induced in cohomology satisfies $s^*_n=\text{id}_{H^*(\ZZ/p^n,\fp)} \otimes \epsilon$, where $\epsilon:H^*(\ZZ/p^n,\fp)\to \fp$ is the augmentation. Therefore,$s_n^*\omega(n)=0$ and $s_n^*c(n)=0$.   
Finally, 
the diagram
\[
\xymatrix{
B_p(n+1) \ar@{->}[r] & B_p(n) \\
\ZZ/{p^{n+1}} \ar@{->}[r] \ar@{->}[u]^{s_{n+1}} & \ZZ/{p^{n}} \ar@{->}[u]^{s_n}
}
\]
commutes and induces a commutative diagram in cohomology
\begin{equation}\label{eq:section square in coh}
\vcenter{
\xymatrix@=1em{
H^*(B_p(n),\fp) \ar[rr]^{\inf_{B_p(n)}^{B_p(n+1)}} \ar[dd]_{s_n^*} && H^*(B_p(n+1),\fp) \ar[dd]^{s_{n+1}^*}\\ \\
H^*(\ZZ/p^n, \fp) \ar[rr]^{\inf_{\ZZ/p^n}^{\ZZ/p^{n+1}}} && H^* (\ZZ/p^{n+1}, \fp).
}
}
\end{equation}
Using this diagram we find that
$$
\alpha_3 c_x(n+1)=s_{n+1}^*\inf_{B_p(n)}^{B_p(n+1)} (c(n)) = \inf_{\ZZ/p^n}^{\ZZ/p^{n+1}} s_n^* (c(n))=0,$$ which gives $\alpha_3c_x(n+1)=0$. Therefore $\alpha_3 =0$, so that $\inf_{B_p(n)}^{B_p(n+1)} c(n) = \alpha_2 \omega_x(n+1)\omega(n+1)$ and the inflated class fully decomposes.

The argument for $P(n)=A_2(n;1)=(\ZZ/2^n)\rtimes \ZZ/2$ is similar. The only class that does not fully decompose in degree~$2$ is $c(n)$. We have
\[
\inf_{A_2(n;1)}^{A_2(n+1;1)}(c(n))
=\alpha_1\,c(n+1)+\alpha_2\,\omega_x(n+1)\,\omega(n+1)+\alpha_3\,\omega_x(n+1)^2.
\]
By considering the inflation--restriction diagram (\ref{eq:inf res square}) for $A_2(n;1)$, we obtain $\alpha_1=0$, and hence the inflation of $c(n)$ is fully decomposable.
\end{proof}
 
\subsection{Proof of \Cref{thm:example}}
 
Next, we generalize the cohomology calculation to all $d\geq 1$, where $d$ is the rank of the normal subgroup in the semidirect products $A_2(n;d)$ and $B_p(n;d,k)$. We start with a general result about group extensions.
 Consider an arbitrary extension  $1\rightarrow N\rightarrow G\rightarrow K\rightarrow 1$ where the action of $K$ on the cohomology of $N$ is trivial, and assume that  in mod $p$ coefficients the cohomology  spectral sequence
 $$
E_2^{p,q}=H^p(K,H^q(N, \fp))\Longrightarrow H^{p+q}(G, \fp) 
 $$ 
 collapses at the $E_2$-page.
\begin{lem}\label{ad}The spectral sequence associated to the pull-back $E$ in the diagram
\begin{equation}\label{ex}
\xymatrix@=1em{
1 \ar@{->}[rr] \ar@{=}[dd]  && N^d \ar@{->}[rr]\ar@{=}[dd] && E \ar@{->}[rr]\ar@{->}[dd]^{\phi} && K \ar@{->}[dd]^{\Delta}\ar@{->}[rr] && 1 \ar@{=}[dd]\\
&&&&&\ar@{=>}[ul]\\
1\ar@{->}[rr] && N^d \ar@{->}[rr] && G^d \ar@{->}[rr] && K^d \ar@{->}[rr] && 1
}
\end{equation}
 collapses at the $E_2$-page, where $\Delta$ is the diagonal map.
\end{lem}
\begin{proof}
 
 

We will compare the spectral sequences associated to the extensions $G^d$ and $E$ whose $E_2$-pages are denoted by $(E')_2^{*,*}$ and $(E'')_2^{*,*}$, respectively. By naturality the vertical maps in (\ref{ex}) induce a map of spectral sequences $(E')_2^{*,*}\rightarrow (E'')_2^{*,*}$. The induced map between the degree $(p,q)$ part of the spectral sequences is 
$$
\phi_{p,q}: H^p(K^d, \fp)\otimes  H^q(N^d, \fp)\rightarrow  H^p(K, \fp)\otimes  H^q(N^d, \fp), 
$$
and it commutes with the differential of the spectral sequence
$$
\xymatrix{
(E')_2^{p,q} \ar@{->}[r]^{\phi_{p,q}} \ar@{->}[d]^{d'_2} & (E'')_2^{p,q} \ar@{->}[d]^{d''_2}\\
(E')_2^{p+2,q-1} \ar@{->}[r]^{\phi_{p+2,q-1}} & (E'')_2^{p+2,q-1}.
}
$$
 The map $\phi_{p,q}$ is surjective since the diagonal map $\Delta$ splits by the projection map $\pi: K^d\rightarrow K$ onto one of the factors.  Hence the composition
$$
H^*(K, \fp)\stackrel{\pi^*}{\rightarrow} H^*(K^d, \fp)\stackrel{\Delta^*}{\rightarrow }H^*(K, \fp)
$$
is the identity map.
We will show that the spectral sequence of $E$ collapses at the $E_2$ page; that is, we show that $d'_2=0$. By the tensor product structure of the terms in the  $E_2$-page, all the differentials are determined by the ones on $(E')_2^{0,q}$. It suffices to show that the edge homomorphism $H^q(G^d, \fp)\rightarrow (E')_2^{0,q}\cong H^q(N^d, \fp)$ is surjective. But this is indeed the case since it is the $d$-fold tensor product of the edge homomorphism $H^q(G, \fp)\rightarrow E_2^{0,q}=H^q(N, \fp)$ which is surjective by the assumption on the spectral sequence of the extension describing  $G$. 
\end{proof}

Since the spectral sequence associated with the extension describing $E$ collapses at the $E_2$-page (by Lemma~\ref{ad}),
the spectral sequence gives an isomorphism of $H^*(K, \fp)$-modules:
$$
H^*(E, \fp)\cong H^*(K, \fp)\otimes H^*(N^d, \fp).
$$
Next, we describe the product structure of the cohomology ring. Consider the diagram 
$$
\xymatrix{
 H^*(K^d,\fp) \ar@{->}[r] \ar@{->}[d]^{\Delta^*} &  H^*(G^d, \fp) \ar@{->}[d]^{\phi^*}\\
H^*(K, \fp) \ar@{->}[r]  & H^*(E, \fp)
}
$$
 induced by (\ref{ex}). The K\"{u}nneth isomorphism \cite[Theorem 2.7.1]{B} yields the following commutative diagram
 $$
\xymatrix{
  \bigotimes^d H^*(K, \fp) \ar@{->}[r]^{\cong} \ar@{->}[rd]^{m} & H^*(K^d, \fp)  \ar@{->}[d]^{\Delta^*} \\
   & H^*(K, \fp)
}
$$
where  $m$ is the multiplication map $a_1\otimes a_2\otimes \cdots \otimes a_d\mapsto a_1a_2\cdots a_d$. Using K\"{u}nneth's theorem for $H^*(G^d, \fp)$, we obtain the commutative diagram of $\bigotimes^d H^*(K, \fp)$-modules
$$
\xymatrix{
\bigotimes^d H^*(K, \fp) \ar@{->}[r] \ar@{->}[d]^{m} & \bigotimes^d H^*(G, \fp) \ar@{->}[d]^{\phi^*}\\
H^*(K,\fp) \ar@{->}[r] & H^*(E, \fp).
}
$$ 
Thus, $\phi^*$ factors as the composition of graded ring homomorphisms
$$
\bigotimes^d H^*(G, \fp)\rightarrow \bigotimes^d_{H^*(K, \fp)} H^*(G, \fp) \stackrel{\varphi}{\longrightarrow} H^*(E, \fp)
$$
where the middle tensor product represents the $d$-fold tensor product of $H^*(G,\fp)$ over $H^*(K,\fp)$:
$$\bigotimes^d_{H^*(K, \fp)} H^*(G, \fp):=\underbrace{H^*(G,\fp) \otimes_{H^*(K,\fp)}\cdots\otimes_{H^*(K,\fp)}H^*(G,\fp)}_{d\text{ times}}.$$
(Recall that $H^*(G,\fp)$ carries an $H^*(K,\fp)$-module structure via the inflation map.)

\Cor{
There is an isomorphism of graded rings
$$
\varphi:\bigotimes^d_{H^*(K, \fp)} H^*(G, \fp)\rightarrow H^*(E, \fp).
$$

}
\Proof{As a graded vector space, $\bigotimes^d_{H^*(K, \fp)} H^*(G, \fp)$ is isomorphic to $H^*(K, \fp)\otimes H^*(N^d, \fp)$ which is isomorphic to $H^*(E, \fp)$  by Lemma \ref{ad}.
} 
 
We apply the above Corollary to $E=A_2(n;d)$ and $E=B_p(n;d,k)$. Let us choose a set of generators $y_1,y_2,\cdots,y_d$ for $(\ZZ/p^n)^d$. Then for any $d\geq 1$ we have
\begin{align*}
H^*(A_2(n;d),\FF_2)&= \FF_2[\omega_i,c_i,\omega_x]/(\omega_i^2+\omega_x\omega_i) \\
H^*(B_p(n;d,k),\FF_p)&= \begin{cases} 
\bigwedge(\omega_x)\otimes \FF_2[\omega_i,c_i,c_x]/(\omega_i^2+\omega_x\omega_i) &\text{if } p=2 \text{ and } k=1 \\
\bigwedge(\omega_x,\omega_i)\otimes \FF_p[c_x,c_i]&\text{if } p \text{ is odd or } p=2 \text{ and } k>1 \end{cases} \end{align*}
where $i=1,2,\cdots,d$. The classes $\omega_i$ and $c_i$ restrict to $\omega_{y_i}$ and $c_{y_i}$ under $\res^E_{(\ZZ/p^n)^d}:H^*(E,\FF_p)\to H^*((\ZZ/p^n)^d,\fp)$, respectively.
We know that under the inflation map $\inf_{\ZZ/p^m}^E:H^*(\ZZ/p^m,\fp)\to H^*(E,\fp)$, every two-dimensional class either maps to zero or fully decomposes into a product of one-dimensional classes (by  \Cref{dec}). One-dimensional classes $\omega_i(n)$ and $\omega_x(n)$ are sent to $\omega_i(n+1)$ and $\omega_x(n+1)$, respectively. 

Therefore in the case where $G_F(p)$ is finitely generated, the cohomology ring $H^*(G_F(p),\FF_p)$ for  superpythagorean fields and $p$-rigid fields is given via
\begin{equation} 
\label{result}
\begin{aligned}
\dlim{}H^*(A_2(n;d),\FF_2) &= \FF_2[\omega_1,\cdots,\omega_d,\omega_x]/(\omega_i^2+\omega_x\omega_i) \\
\dlim{}H^*(B_p(n;d,k),\FF_p)&= \begin{cases} 
\bigwedge(\omega_x)\otimes \FF_2[\omega_1,\cdots,\omega_d]/(\omega_i^2+\omega_x\omega_i) &\text{if } p=2 \text{ and } k=1 \\
\bigwedge(\omega_x,\omega_1,\cdots,\omega_d)&\text{if } p \text{ is odd or } p=2 \text{ and } k>1. \end{cases}
\end{aligned} 
\end{equation}
See also \cite{Wadsworth} for a different approach to this calculation.

To conclude the proof of Theorem \ref{thm:example}, first consider the case where $F$ is a superpythagorean field. Let $\gamma \in H^2(P(n),\mathbb{F}_2)$ be given. Then from the formula $$G_F(2) = \prod_\mathcal{I} \mathbb{Z}_2\rtimes C_2$$ and $P(n)=\prod_\mathcal{I} C_{2^n}\rtimes C_2$, we see that $\gamma$ is in fact the inflation of some $\gamma' \in H^*\left(\left(\prod_\mathcal{J} C_{2^n}\right)\rtimes C_2,\mathbb{F}_2\right)$ where $\mathcal{J}$ is finite. From our results above and the transitivity of inflation maps, we can conclude that the inflation of $\gamma$ to $H^2(P(n+1),\mathbb{F}_2)$ fully decomposes. 

The proof when $F$ is a $p$-rigid field is similar.

This concludes the proof of Theorem \ref{thm:example}, a special case of Theorem \ref{thm:main}, which relies on group cohomological methods independent of the Bloch--Kato conjecture.

\section{Groups whose cohomology ring is generated in degree one}\label{sec:BKpropgroups}
Having studied some decomposability results for Galois cohomology, we now characterize those finite groups $G$ with trivial action on $\fp$ for which the cohomology ring $H^*\left( G, \fp\right)$ is generated in degree one.  These, and similar groups discussed below, were first introduced in an earlier version of this paper and were investigated by other authors; see for example \cites{Q14,QW20}.

We first introduce some notation that will be necessary to state the main result of this section. Let $G$ be any finite group, and let $p$ be a prime number. We define $O_{p'}(G)$ to be the subgroup of $G$ generated by those elements of order coprime to $p$; that is,
$$O_{p'}(G) = \langle \{g \in G: \gcd(|g|,p)=1\}\rangle.$$ 
Since isomorphisms preserve orders of elements, it follows that for any automorphism $\psi:G \to G$ we have $\psi(O_{p'}(G)) = O_{p'}(G)$; in particular we have that $O_{p'}(G)$ is a normal subgroup of $G$. Furthermore, if $H$ is a $p$-group and $\Delta:G \to H$ is a group homomorphism, then $O_{p'}(G) \subseteq \ker(\Delta)$. Since $G/O_{p'}(G)$ is a $p$-group, this implies that $O_{p'}(G)$ is the smallest normal subgroup of $G$ whose quotient is a $p$-group. 

In the statement below, for a prime $p$, we will make reference to the class of elementary $p$-abelian groups. These are abelian groups with the property that every element has order dividing $p$. 

\begin{thm}\label{thm:cohomology.generated.degree.one}
    Suppose that $G$ is a finite group,  that $p$ is a prime number dividing $|G|$, and that $G$ acts trivially on $\fp$. The following statements are equivalent:
    \begin{enumerate}
        \item $H^*(G, \fp)$ is generated in degree one; 
        \item $p=2$, and $G/O_{2'}(G)$ is a nontrivial elementary $2$-abelian group, and $|O_{2'}(G)|$  is odd;
        \item $p=2$, and $G \cong N \rtimes \prod_{i=1}^n C_2$ for some $n \geq 1$ and $N$ a finite group with $|N|$ odd; and
        \item $p=2$, and the Sylow $2$-subgroup of $G$ is a nontrivial elementary $2$-abelian group that admits a normal complement.  
    \end{enumerate}
\end{thm}

We remark that for any group that satisfies the equivalent conditions of the theorem, $O_{2'}(G)$ is solvable by the Feit--Thompson theorem \cite[\S 1.1, p. 775, Theorem]{Feit-Thompson}, and therefore $G$ must be solvable as well by (3). (Solvable groups are closed under group extensions.) For an example of a family of groups which satisfy this condition, one can consider the dihedral group of order $2n$ where $n$ is an odd number. The case where $n=3$ tells us that $S_3$ is one such group. Finally, we note that if $p$ does not divide the order of $G$, then $H^*(G, \fp)$ is concentrated in degree $0$. 

\begin{proof}
Let $\phi_p(G)$ be the intersection of all normal subgroups of $G$ of index $p$.  Equivalently, this is the intersection of the kernels of all homomorphisms from $G$ to $\fp$; that is
$$\phi_p(G) = \bigcap_{\theta\in\Hom(G,\fp)} \ker(\theta).$$ 
We observe that $O_{p'}(G) \subseteq \phi_p(G) \subseteq G$.  

$(1) \Rightarrow (2)$: From (1) we know that $H^*(G,\fp)$ is generated by degree-one elements. We show that each element of $\phi_p(G)$ has order coprime to $p$. Suppose instead that there exists some $g \in \phi_p(G)$ 
so that $|g|$ is not relatively prime to $p$. Hence a suitable power of $g$ --- let us call it $h$ --- has order precisely $p$.   From \cite[Theorem 1, p.~885]{swan}, we conclude that the restriction map $$\res: H^*(G,\fp) \to H^*(\langle h \rangle,\fp)$$ is nontrivial for an infinite number of degrees.  Since $h \in \phi_p(G)$, for any $\theta \in H^1(G, \fp)$ we have $\res(\theta) = 1$. Because $H^*(G,\fp)$ is generated by degree one elements, we conclude that the restriction map is trivial in all positive degrees, a contradiction. Hence every $g \in \phi_p(G)$ has order relatively prime to $p$, and so we have $\phi_p(G) \subseteq O_{p'}(G)$. From our previous observation it follows that $\phi_p(G)=O_{p'}(G)$. 

Now let us show that $G/O_{p'}(G) = G/\phi_p(G)$ is a nontrivial, elementary $p$-abelian group.  For nontriviality, observe that since $p \mid |G|$ we know from Cauchy's theorem that $G$ contains an element of order $p$.  Hence $\phi_p(G)$ --- which we showed above consists only of elements of order coprime to $p$ --- cannot be all of $G$, and so $G/\phi_p(G)$ is nontrivial. To see that $G/\phi_p(G)$ is elementary $p$-abelian, let $\mathcal{H}$ be the collection of normal subgroups of index $p$ within $G$, i.e., kernels of homomorphisms $\theta$ as defined in the intersection above.


If we define $$E = \prod_{H \in \mathcal{H}} G/H \simeq \prod_{H \in \mathcal{H}} \mathbb{F}_p,$$ then the natural map $G \to E$ has kernel $\phi_p(G)$.

We now proceed to show that $p$ is even, which we do by contradiction.  Assume then that $p$ is odd.  From above we know that $|O_{p'}(G)|$ is relatively prime to $p$, so we deduce that $H^*(O_{p'}(G),\mathbb{F}_p) = \mathbb{F}_p$ is concentrated in degree $0$.  Hence considering the Lyndon--Hochschild--Serre spectral sequence associated with $$1 \to O_{p'}(G) \to G \to G/O_{p'}(G) \to 1,$$ we see that $\text{inf}:H^*(G/O_{p'}(G),\mathbb{F}_p) \to H^*(G,\mathbb{F}_p)$ is an isomorphism. By \cite[Corollary 3.5.7(b)]{B} we have that $H^*(G/O_{p'}(G),\fp)$ contains an element in degree $2$ which is not contained in the subalgebra of $H^*(G/O_{p'}(G),\fp)$ generated in degree $1$. Hence the inflation of this element to $H^*(G,\fp)$ is also not in the subalgebra of $H^*(G,\fp)$ generated by degree $1$. Again we reach a contradiction, and so it follows that $p=2$. 

Because we assume that $H^*(G,\mathbb{F}_2)$ is generated by degree $1$ elements, the same argument as before shows that no element $g \in O_{2'}(G)$ has even order (lest some power of that element have order $2$, in which case we arrive at a contradiction).  Therefore $O_{2'}(G)$ is a subgroup of $G$ of odd order. This shows (2). 

$(2) \Rightarrow (1)$ Since $|O_{2'}(G)|$ is odd, we deduce that $H^*(O_{2'}(G),\mathbb{F}_2) = \mathbb{F}_2$ is concentrated in degree $0$.  Hence considering the Lyndon--Hochschild--Serre spectral sequence associated with $$1 \to O_{2'}(G) \to G \to \prod_{k=1}^n \mathbb{F}_2 \to 1$$ for some $n \in \mathbb{N}$, we see that $\text{inf}:H^*(G/O_{2'}(G),\mathbb{F}_2) \to H^*(G,\mathbb{F}_2)$ is an isomorphism.  Since $G/O_{2'}(G)$ is a nontrivial elementary $2$-abelian group, we therefore conclude that $H^*(G,\mathbb{F}_2)$ is a polynomial algebra over $\mathbb{F}_2$ generated by degree $1$; see \cite{B2}.

The Schur--Zassenhaus theorem (see \cite[Ch.~5, pg.~37]{Wehr}) implies $(2) \Rightarrow (3)$, and $(3) \Rightarrow (4)$ is immediate. For $(4) \Rightarrow (2)$, suppose that $S$ is the Sylow $2$-subgroup of $G$, that $S$ is elementary $2$-abelian, and that $S$ admits a normal complement $N$. Then $|N|$ is odd. In particular this means that $N \subseteq O_{2'}(G)$.  We know from the discussion before the statement of Theorem \ref{thm:cohomology.generated.degree.one} that $O_{2'}(G)$ is the smallest subgroup of $G$ such that $G/O_{2'}(G)$ is a $2$-group; on the other hand, we know that $G/N \simeq S$ since we have assumed that $N$ is a normal complement to $S$.  Hence we conclude that $N = O_{2'}(G)$. Therefore we have $|O_{2'}(G)|$ is odd, and $G/O_{2'}(G) \simeq S$ is elementary $2$-abelian.
\end{proof}

\begin{rem} See \cite[Ch.~5, Ex.~5]{Wehr} for other consequences of the property $p \nmid O_{p'}(G)$ of the group $G$.
\end{rem}

A direct consequence of this theorem is a characterization of finite $p$-groups so that the cohomology ring is generated in degree one. This answers a question posed in \cite[Remark 2.10]{Q22}. 

\begin{cor}
 \label{Claudio}
    Let $G$ be a nontrivial finite $p$-group. The cohomology ring $H^*(G, \fp)$ is generated by degree-one elements if and only if $p=2$ and $G$ is an elementary $2$-abelian group.
\end{cor}

\addresseshere
\newpage
\appendix
\thispagestyle{empty}
\section{Computations for \texorpdfstring{$Q_8$}{Q8} \texorpdfstring{\\}{}by David Benson}\label{appendix}

Let $Q_8$ be the quaternion group of order eight,
\[ Q_8=\langle g,h\mid g^4=1,\ g^2=h^2,\ gh=h^3g\rangle. \]
It is well known that \cite{Swan:1960}
\[ H^*(Q_8, \bF_2) = \bF_2[z,y,v]/(y^2+yz+z^2,y^2z+yz^2), \]
where $|z|=|y|=1$, $|v|=4$. 
Taking cup product with the class $v$ gives an isomorphism on cohomology groups:
\[ \cup \;v \colon H^n(Q_8, \bF_2) \overset{\cong}{\longrightarrow} H^{n+4}(Q_8, \bF_2) \; \; \forall \; n > 0.\]
Thus, the cohomology is periodic with period four
and looks as
follows:
\[ \begin{array}{|c|c|c|c|c|c|c|c|} \hline
\deg & 0 & 1 & 2 & 3 & 4 & 5 & 6 \\ \hline
\text{basis}&1&y,z&y^2,yz&y^2z&v&vy,vz & \dots\\ \hline
\end{array} \]

The purpose of this appendix is to prove the following.

\begin{proposition}\label{pr:Q8}
Let
\[ H = \langle g,h\mid g^8=1,\ g^4=h^4,\ gh=h^3g \rangle, \]
a group of order 32,
and let $H \to Q_8$ be the surjective homomorphism taking $g$ and $h$ to
the elements with the same name in $Q_8$. Then $\Inf_{Q_8}^H(v)=0$.

Let $H\times\bZ/2\to Q_8$ be the homomorphism with first factor as
above, and where $\bZ/2$ is taken isomorphically to the centre of
$Q_8$. Then $\Inf_{Q_8}^{H\times\bZ/2}(v)$ is non-zero, but is in the subalgebra
generated by degree one elements.
\end{proposition}

For general background on cohomology of groups, we refer to Chapter~12
of Cartan and Eilenberg~\cite{Cartan/Eilenberg:1956a} and to Adem and
Milgram~\cite{Adem/Milgram:1994a}. 
For the techniques used in the proof, we refer
to Benson~\cite{Benson:1991b}, where an extensive discussion of
Steenrod operations in group cohomology occurs in Section~4.4; see
also May~\cite{May:1970a}. We make
use of Kudo's transgression 
theorem~\cite{Kudo:1956a}, which states that Steenrod operations
commute with transgressions, to determine some of the differentials in the spectral
sequence of a group extension as in~\cite[\S4.1]{Adem/Milgram:1994a},
\cite[\S4.8]{Benson:1991b}. The properties of the Steenrod 
operations used for computation are listed
in~\cite[Theorem~4.4.4]{Benson:1991b}; in particular, the Cartan
formula says that it suffices to know the action of Steenrod
operations on generators. 

We remark that extensive computations using 
{\sf Magma}~\cite{Bosma/Cannon/Playoust:1997a} show that 
$H$ is the unique group of order at most 32 with the
above property. It seems likely that if $G\to Q_8$ is a surjection with the property
that $\Inf_{Q_8}^G(v)=0$ then $G\to Q_8$ factors as $G\to H\to Q_8$ with $G\to H$ a
surjection and $H\to Q_8$ as above. Furthermore, it seems likely that if 
$\Inf_{Q_8}^G(v)$ is non-zero, but is in the subalgebra generated by degree one elements 
then $G\to Q_8$ factors as $G \to H\times \bZ/2 \to
Q_8$ with $G\to H\times\bZ/2$ a surjection and $H\times\bZ/2\to Q_8$
as above.

Throughout, we have abused notation by using the same names for
elements of different groups, and of different cohomology rings. This
should not cause confusion.

\begin{proof}[Proof of Proposition~\ref{pr:Q8}]
Let 
\[ V_4=\langle g,h\mid g^2=1,\ h^2=1,\ gh=hg\rangle\cong\bZ/2\times\bZ/2, \]
with cohomology 
\[ H^*(V_4,\bF_2)=\bF_2[z,y] \]
where $|z|=|y|=1$. Let
\[ H_1=\langle g,h\mid g^4=1,\ h^2=1,\ gh=hg\rangle\cong\bZ/2\times\bZ/4, \] 
with cohomology
\[ H^*(H_1,\bF_2)=\bF_2[z,y,x]/(z^2) \]
where $|z|=|y|=1$, $|x|=2$, $\Sq^1(z)=0$, $\beta_2(z)=x$,
$\Sq^1(y)=y^2$, $\Sq^1(x)=0$. 
Here, $\beta_2$ denotes the second Bockstein~\cite[\S4.3]{Benson:1991b}.
The notation is chosen so that $z$ and $x$ are inflated from $H_1/\langle h\rangle$ and $y$ is 
inflated from $H_1/\langle g\rangle$.

Let
\[ H_2=\langle g,h\mid g^4=1,\ h^4=1,\ gh=h^3g\rangle \cong \bZ/4\rtimes\bZ/4. \]
This is SmallGroup(16,4), Hall/Senior~\cite{Hall/Senior:1964a} group of order $16$ 
number $10$, named $\Gamma_2c_2$. Then $H_2$ 
surjects onto $H_1$ by sending group elements to those with the same name. 
Examining squares and commutators, the class of the extension 
is $y^2+yz$. This is because $h$ lifts to an element of order $4$ and 
$[g,h]$ is now $h^2$ (in detail, it follows by restricting to
$\langle g\rangle$ and $\langle h\rangle$ that the extension class is either $y^2$ or
$y^2+yz$; but the extension corresponding to $y^2$ is the abelian
group $\bZ/4\times\bZ/4$).  In the spectral sequence 
\[ H^*(H_1,H^*(\bZ/2,\bF_2))\Rightarrow H^*(H_2,\bF_2) \]
let $H^*(\bZ/2,\bF_2)=\bF_2[\zeta]$ with $|\zeta|=1$. Then $d_2(\zeta)=y^2+yz$, which is not
a zero divisor in $H^*(H_1,\bF_2)$. So $E_3=\bF_2[z,y,x,\zeta^2]/(z^2,y^2+yz)$
and by Kudo's transgression theorem~\cite{Kudo:1956a},
$d_3(\zeta^2)=\Sq^1(y^2+yz)=y^2z$ which is already zero in the $E_3$ page.
Let $w$ be a representative of $\zeta^2$. Then there are no relations to ungrade, so
\[ H^*(H_2,\bF_2)=\bF_2[z,y,x,w]/(z^2,y(y+z)) \]
with $|z|=|y|=1$, $|x|=|w|=2$.
The elements $z$, $y$ and $x$ are inflations from $H_1$ of elements with the same names.
The element $w$ has non-zero restriction to the central element $\langle h^2\rangle$.

There is a surjective map $H_2\to Q_8$ with
kernel of order two generated by $g^2h^2$. This gives us a diagram
of groups and homomorphisms:
\[ \xymatrix{1\ar[r]&\bZ/2\ar[r]\ar@{=}[d]&\bZ/4\rtimes\bZ/4\ar[r]\ar[d]&\bZ/2\times\bZ/4 \ar[r]\ar[d]&1\\
1\ar[r]&\bZ/2\ar[r]&Q_8\ar[r]&\bZ/2\times\bZ/2\ar[r]&1,} \]
i.e.,
\[ \xymatrix{1\ar[r]&\bZ/2\ar[r]\ar@{=}[d]&H_2\ar[r]\ar[d]&H_1\ar[r]\ar[d]&1\\
1\ar[r]&\bZ/2\ar[r]&Q_8\ar[r]&V_4\ar[r]&1.} \]
All the maps in the right hand square take elements $g$ and $h$ to the elements with the 
same names. This gives us a map of spectral sequences
\[ \xymatrix{H^*(V_4,H^*(\bZ/2,\bF_2)) \ar@{=>}[r]\ar[d]&H^*(Q_8,\bF_2)\ar[d]\\
H^*(H_1,H^*(\bZ/2,\bF_2))\ar@{=>}[r]&H^*(H_2,\bF_2).} \]
In the $E_2$ pages, we let $H^*(\bZ/2,\bF_2)=\bF_2[\zeta]$ with $|\zeta|=1$. Thus the map of 
$E_2$ pages is as follows:
\[ \bF_2[z,y,\zeta] \to \bF_2[z,y,x,\zeta]/(z^2) \]
where, again, elements map to elements with the same name, by abuse of notation.
We have $d_2(\zeta)=y^2+yz+z^2$ in $H^*(V_4,\bF_2)$, which inflates to $d_2(\zeta)=y^2+yz$ in
$H^*(H_1,\bF_2)$ since $z^2=0$ there. In both cases, this is a non zero divisor, so at the
level of $E_3$ pages we have
\[ \bF_2[z,y,\zeta^2]/(y^2+yz+z^2) \to \bF_2[z,y,x,\zeta^2]/(z^2,y^2+yz). \]
Now $\Sq^1(y^2+yz+z^2)=y^2z+yz^2$. This is non-zero on the left but zero on the right.
So by Kudo's transgression theorem we have $d_3(\zeta^2)=y^2z+yz^2$ on the left and
$d_3(\zeta^2)=0$ on the right. So the $E_4=E_5$ page on the left is
\[ \bF_2[z,y,\zeta^4]/(y^2+yz+z^2,y^2z+yz^2) \to \bF_2[z,y,x,\zeta^2]/(z^2,y^2+yz). \]
Since
\[ \Sq^2(y^2z+yz^2)=y^4z+yz^4= (y^2+yz+z^2)(y^2z+yz^2) \]
we have $d_5(\zeta^4)=0$, so $E_5=E_\infty$, and both spectral sequences have stopped.

Let $w$ be a representative of $\zeta^2$ in $H^*(H_2,\bF_2)$, and let $v$ be a 
representative of $\zeta^4$ in $H^*(Q_8,\bF_2)$.  Thus
\begin{align*} 
H^*(H_2,\bF_2)&=\bF_2[z,y,x,w]/(z^2,y^2+yz),\\
H^*(Q_8,\bF_2)&=\bF_2[z,y,v]/(y^2+yz+z^2,y^2z+yz^2). 
\end{align*}
with $|z|=|y|=1$, $|x|=|w|=2$, $|v|=4$.

Using the map of spectral 
sequences, the choice can be made so that $v\in H^4(Q_8)$ inflates to
$w^2\in H^4(H_2)$, as follows.
The restriction of $\Sq^1(w)$ to $\bZ/2\times\bZ/4$ is zero, so
$\Sq^1(w)$ is a linear combination of $zw$ and $zx$. This implies that
$\Sq^2(w^2)=(\Sq^1(w))^2=0$. 
The map of spectral sequences shows that the inflation of $v$ from
$Q_8$ to $H_2$ is some representative of the element $\zeta^4$. 
The indeterminacy in $\zeta^4$ is the linear
span of $\zeta^2y^2$, $\zeta^2x$, $xy^2$ and $x^2$. However, the
inflation of $\Sq^2(v)$ has to be zero, while $\Sq^2(w^2)=0$,
$\Sq^2(wy^2)=w^2y^2$, $\Sq^2(xy^2)=x^2y^2$, $\Sq^2(x^2)=0$. Thus the
inflation of $v$ from $Q_8$ to $H_2$ is either $w^2$ or $w^2+x^2$.
But the indeterminacy in $w$ allows us to replace it by $w+x$ if
necessary so that the inflation of $v$ is $w^2$.

Now let $H$ be the extension
\[ 1\to \bZ/2\to H \to H_2 \to 1 \]
corresponding to $w\in H^2(H_2,\bF_2)$. Then $w$ inflates to zero in $H$, and
hence $v$ inflates to $w^2$ in $H^*(H_2)$ and then to zero in $H^*(H)$.  The group $H$ has presentation
\[ H=\langle g,h\mid g^8=1,\ g^4=h^4,\ gh=h^3g\rangle. \] 
This is SmallGroup(32,15), Hall/Senior~\cite{Hall/Senior:1964a} group 
of order 32 number 32, named $\Gamma_3f$. This completes the proof of
the first part of Proposition~\ref{pr:Q8}.

For the second part, we note that the map $H_2\times\bZ/2 \to Q_8$
factors as $H_2\times \bZ/2 \to Q_8\times \bZ/2 \to Q_8$, where the
map $Q_8\times \bZ/2\to Q_8$ takes $g$ to $g$, $h$ to $h$, and the
generator of $\bZ/2$ to the central element $g^2=h^2$. Since elements
of $H^*(Q_8)$ of degree less than four are in the subalgebra generated
by degree one elements, the inflation of $v$ to $Q_8\times \bZ/2$ is
$v\otimes 1$ plus elements in the subalgebra generated by degree one
elements. Then inflating further to $H_2\times\bZ/2$, the term
$v\otimes 1$ inflates to zero, by the first part of the Proposition.
To show that the inflation is non-zero, we notice that if we inflate
from $Q_8$ to $H\times\bZ/2$ and then restrict to $\bZ/2$, we get the same answer
as if we just restrict directly from $Q_8$ to the central
$\bZ/2$, regarded as a subgroup of $H\times\bZ/2$. 
This restriction sends $v$ to the fourth power of the
generator, which is non-zero.
\end{proof}

\begin{remark}
Note that $H$ is a central extension of $\bZ/4\rtimes\bZ/4$, and is a 
quotient by a central element $g^4h^4$ of order two of a
semidirect product 
\[ \bZ/8\rtimes\bZ/8 = \langle g,h\mid g^8=1,\ h^8=1,\
  gh=h^3g\rangle. \]
It follows that the inflation of $v$ from $Q_8$ to the latter
semidirect product is zero.
\end{remark}

\newcommand{\noopsort}[1]{}
\providecommand{\bysame}{\leavevmode\hbox to3em{\hrulefill}\thinspace}
\providecommand{\MR}{\relax\ifhmode\unskip\space\fi MR }
\providecommand{\MRhref}[2]{%
  \href{http://www.ams.org/mathscinet-getitem?mr=#1}{#2}
}
\providecommand{\href}[2]{#2}

\end{document}